\documentclass[preprint,12pt]{elsarticle}
\usepackage[utf8]{inputenc}
\usepackage{times}
\usepackage{amsthm,amsmath,amsfonts,latexsym,amssymb}
\usepackage{graphics,fancyhdr,graphicx}
\usepackage{multirow}
\usepackage{float}
\usepackage{fleqn}
\usepackage{subfigure}
\usepackage{verbatim}
\usepackage{rotating}
\usepackage{makecell}
\usepackage{bm}
\usepackage{color}
\usepackage{bm}

\usepackage{colortbl}
\usepackage[table]{xcolor}

\usepackage{geometry}
\biboptions{numbers,sort&compress}
\definecolor{listinggray}{gray}{0.9}
\definecolor{lbcolor}{rgb}{0.9,0.9,0.9}
\definecolor{Darkgreen}{RGB}{0,100,0}

\begin{document}
\abovedisplayskip=6.0pt
\belowdisplayskip=6.0pt
\begin{frontmatter}
\title{Generalized Pole-Residue Method for Dynamic Analysis of Nonlinear Systems based on Volterra Series}

\author[dlut]{Qianying Cao}
\ead{caoqianying@dlut.edu.cn}
\author[ouc]{Anteng Chang}
\ead{changanteng@ouc.edu.cn}
\author[ouc]{Junfeng Du}
\ead{dujunfeng@ouc.edu.cn}
\author[dlut]{Lin Lu\corref{cor1}}
\ead{lulin@dlut.edu.cn}

\cortext[cor1]{Corresponding author.}

\address[dlut]{State Key Laboratory of Coastal and Offshore Engineering, Dalian University of Technology, Dalian 116024, PR China}
\address[ouc]{College of Engineering, Ocean University of China, Qingdao 266100, PR China}

\begin{abstract}
Dynamic systems characterized by second-order nonlinear ordinary differential equations appear in many fields of physics and engineering. To solve these kinds of problems, time-consuming step-by-step numerical integration methods and convolution methods based on Volterra series in the time domain have been widely used. In contrast, this work develops an efficient generalized pole-residue method based on the Volterra series performed in the Laplace domain. The proposed method involves two steps: (1) the Volterra kernels are decoupled in terms of Laguerre polynomials, and (2) the partial response related to a single Laguerre polynomial is obtained analytically in terms of the pole-residue method. Compared to the traditional pole-residue method for a linear system, one of the novelties of the pole-residue method in this paper is how to deal with the higher-order poles and their corresponding coefficients. Because the proposed method derives an explicit, continuous response function of time, it is much more efficient than traditional numerical methods. Unlike the traditional Laplace domain method, the proposed method is applicable to arbitrary irregular excitations. Because the natural response, forced response and cross response are naturally obtained in the solution procedure, meaningful mathematical and physical insights are gained. In numerical studies, systems with a known equation of motion and an unknown equation of motion are investigated. For each system, regular excitations and complex irregular excitations with different parameters are studied. Numerical studies validate the good accuracy and high efficiency of the proposed method by comparing it with the fourth-order Runge--Kutta method.
\end{abstract}


\end{frontmatter}

\section{Introduction}
Most real dynamic systems, as encountered in mechanical and civil engineering, are inherently nonlinear and include geometric nonlinearities, nonlinear constitutive relations in material or nonlinear resistances, etc.~\cite{cheng2017volterra}. Nonlinear problems are attracting increasing attention from engineers and scientists. This work focuses on solving nonlinear system vibration problems, \textit{i.e.}, computing transient responses of nonlinear oscillators under arbitrary irregular excitations based on a combination of a pole-residue operation and Volterra series. Because Volterra series are single-valued, the scope of the present study is restricted to nonlinear behaviours without bifurcations~\cite{worden1997harmonic}.

To analyse nonlinear vibration problems, researchers have performed extensive studies and developed various mathematical methods. Popular methods include step-by-step numerical integration methods in the time domain, such as the Runge--Kutta method. This kind of method not only requires a small time-step resolution for obtaining high-precision solutions but also is prone to numerical instability~\cite{meirovitch1997principles,iserles2009first}. For a long response with small time steps, the time domain methods are very costly in computational time.
Volterra series is another widely used method, which is the extension of the Duhamel integral for linear systems~\cite{volterra1959theory,schetzen1980volterra}. Volterra series can reproduce many nonlinear phenomena, but they are very complex due to higher-dimensional convolution integrals~\cite{cheng2017volterra}. Since 1980's, significant progress has been made in the general area of the Volterra series. The reader is referred to Ref.~\cite{cheng2017volterra} for a quite thorough literature review on the relevant topics. After 2017, most papers focus on Volterra series identification. De Paula and Marques~\cite{de2019multi} proposed a method for the identification of Volterra kernels, which was based on time-delay neural networks. Son and Kim~\cite{son2020parametric} presented a method for a direct estimation of the Volterra kernel coefficients. Dalla Libera et al.~\cite{dalla2021kernel} introduced two new kernels for Volterra series identification. Peng et al.~\cite{peng2021nonlinear} used the measured response to identify the kernel function and performed the nonlinear structural damage detection. Only a few papers concentrated on simplifying the computation of convolution integrals. Traditional methods for computing convolution integrals involved in the Volterra series have been performed in three distinct domains: time, frequency and Laplace. The time domain method based on Volterra series refers to discrete time convolution methods, which also suffer computational cost problems~\cite{liu2012new,son2020parametric}. Both the frequency domain method and the Laplace domain method based on the Volterra series consist of three steps: (1) Volterra series are transformed into an algebraic equation in the frequency domain or Laplace domain; (2) the algebraic equation is solved by purely algebraic manipulations; and (3) the solution in Step (2) is transformed back to the time domain. Many researchers have used the frequency domain method to compute the responses of nonlinear systems. Billings et al.~\cite{billings1989spectral} developed a new method for identifying the generalized frequency response function (GFRF) of nonlinear systems and then predicted the nonlinear response based on these GFRFs. Carassale et al.~\cite{carassale2014nonlinear} introduced a frequency domain approach for nonlinear bridge aerodynamics and aeroelasticity. Ho et al.~\cite{ho2014frequency} computed an output frequency domain function of a nonlinear damped duffing system modelled by a Volterra series under a sinusoidal input. Kim et al.~\cite{kim2014time,kim2015finite} identified the higher order frequency response functions by using the nonlinear autoregressive with exogenous input technique and the harmonic probing method. This type of frequency domain method is much more efficient than the time domain method due to the fast Fourier transform algorithm. However, the frequency domain method not only is limited by frequency resolutions but also suffers from leakage problems due to the use of discrete Fourier transforms. In addition, the frequency domain method calculates only a steady-state response. A natural response generated by initial conditions and a cross response caused by interactions between a system and an excitation are ignored. In contrast, the Laplace domain method can calculate all response components because initial conditions are considered in the computational procedure. However, it has been restricted to analytical operations for simple excitations, such as sinusoidal excitations and exponential excitations~\cite{kreyszig2010advanced}.

The proposed method falls into the category of the Volterra series method computed in the Laplace domain. Unlike the traditional Laplace domain method, the proposed method is applicable to arbitrary irregular excitations. Because the proposed method follows a similar path as a pole-residue method for linear systems~\cite{hu2016pole}, the proposed method to solve nonlinear system vibration problems is called the generalized pole-residue method. The main concept of the pole-residue method developed by Hu et al.~\cite{hu2016pole} was that the poles and residues of the response could be easily obtained from those of the input and system transfer function to obtain the closed-form response solution of linear systems. This method included three steps: (1) writing the system transfer function into pole-residue form; (2) writing the excitation into pole-residue form by the Prony-SS method; (3) computing the poles and residues of the response by an algebraic operation based on those from system and excitation. Compared to Hu et al.~\cite{hu2016pole}, which was regarded as an efficient tool to compute responses of linear systems, the generalized pole-residue method in this paper is introduced to compute responses of nonlinear systems. The proposed method involves two steps: (1) the Volterra kernels are decoupled in terms of Laguerre polynomials, and (2) the partial response related to a single Laguerre polynomial is obtained analytically in terms of the pole-residue method. Compared to the traditional pole-residue method for a linear system, one of the novelties of the generalized pole-residue method is how to deal with the higher-order poles and their corresponding coefficients.
Similar to the Taylor series, the Volterra series representation is an infinite series, and convergence conditions are needed to assure that the representation is meaningful. Because the proposed method is based on the Volterra series, only the system with convergent Volterra series representation can be treated by the proposed method.

The paper is organized as follows. In Section 2, the nonlinear response is modelled by a Volterra series, and Volterra kernel functions are decoupled by Laguerre polynomials. Then, the pole-residue method for computing explicit responses is developed in Section 3. Numerical studies and discussions are given in Section 4. Finally, the conclusions are drawn in Section 5.

\section{Response calculation based on Volterra series}
A nonlinear oscillator, whose governing equation of motion is given by
\begin{equation}\label{eq:equationofmotion}
  m\ddot{y}(t)+c\dot{y}(t)+ky+z(t,y,\dot{y})=f(t)
\end{equation}
where $z(t,y,\dot{y})$ represents an arbitrary nonlinear term; $m$, $c$, and $k$ are the mass, damping and linear stiffness, respectively; $y(t)$, $\dot{y}(t)$ and $\ddot{y}(t)$ are the displacement, velocity and acceleration, respectively; and $f(t)$ is the time-dependent excitation.

If the energy of excitation $f(t)$ is limited, the nonlinear response under zero initial conditions (\textit{i.e.}, zero displacement and zero velocity) can be represented by the Volterra series~\cite{boyd1984analytical,boyd1985fading,cheng2017volterra,rugh1981nonlinear}:
\begin{equation}\label{y1}
  y(t)=\sum_{n=1}^Ny_n(t)
\end{equation}
where $N$ is the order of Volterra series and
\begin{equation}\label{yn}
  y_n(t)=\int_{-\infty}^{\infty}\ldots\int_{-\infty}^{\infty}h_{n}(\tau_{1},\ldots,\tau_{n})\prod_{i=1}^nf(t-\tau_i)d\tau_1\ldots d\tau_n
\end{equation}
In Eq.~\ref{yn}, $h_1(\tau)$ is called the first-order Volterra kernel function, which represents the linear behaviour of the system; $h_{n}(\tau_{1},\ldots,\tau_{n})$ for $n>1$ are the higher-order Volterra kernel functions, which describe the nonlinear behaviour of the system.
The complete formulation of $y(t)$ includes infinite series where the labour of calculating the $n^{th}$ term increases quickly with the growth of $n$. Fortunately, the response accuracy may be ensured by the first several order Volterra series. This is proved here in numerical studies.

The commonly known Laguerre polynomials are represented as~\cite{israelsen2014generalized,son2020parametric}:
\begin{equation}\label{laguerre}
  l_{p_i}(t)=\sqrt{2a_i}\sum_{k=0}^{p_i}\frac{(-1)^{k}p_i!}{k![(p_i-k)!]^{2}}(2a_it)^{p_i-k}e^{-a_it}
\end{equation}
where $p_i$ is the order of the Laguerre polynomials and $a_i$ is the damping rate. The Laguerre polynomials satisfy the orthogonal relationship expressed as:
\begin{equation}\label{eq:orthogonal}
  \int_0^{\infty}l_{p_i}(t)l_{p_j}(t)dt=\delta=\left\{\begin{array}{c}
                                              0~ (p_i\neq p_j)\\
                                              1~(p_i = p_j)
                                            \end{array}\right.
\end{equation}

By using Laguerre polynomials, the Volterra kernel function $h_{n}(t_{1},\ldots,t_{n})$ in Eq.~\ref{yn} can be decoupled as follows~\cite{israelsen2014generalized,son2020parametric}:
\begin{equation}\label{ht_exp}
 h_{n}(t_{1},\ldots,t_{n})=\sum_{p_1=0}^{R_1}\ldots\sum_{p_n=0}^{R_n}c_{p_1\ldots p_n}l_{p_1}(t_1)\ldots l_{p_n}(t_n)
\end{equation}
where the coefficient is computed resorting to the orthogonal relationship in Eq.~\ref{eq:orthogonal}:
\begin{equation}\label{eq:coefficient}
  c_{p_1\ldots p_n}=\int_0^{\infty}\ldots\int_0^{\infty}l_{p_1}(t_1)\ldots l_{p_n}(t_n) h_{n}(t_{1},\ldots,t_{n})dt_1\ldots dt_n
\end{equation}
Substituting Eq.~\ref{ht_exp} into Eq.~\ref{yn} yields
\begin{equation}\label{y2}
  y_n(t)=\sum_{p_1=0}^{R_1}\ldots\sum_{p_n=0}^{R_n}c_{p_1\ldots p_n}\int_{-\infty}^{\infty}\ldots\int_{-\infty}^{\infty}l_{p_1}(\tau_1)\ldots l_{p_n}(\tau_n)\prod_{i=1}^nf(t-\tau_i)d\tau_i
\end{equation}
The above operation that uses the Laguerre polynomials to decouple Volterra higher order kernel functions has been well-developed. The reader is referred to Refs.~\cite{schetzen1980volterra,israelsen2014generalized} for details about the adopted technique. After decoupling Volterra higher order kernel functions in time, one can regroup Eq.~\ref{y2} into:
\begin{equation}\label{y3}
  y_n(t)=\sum_{p_1=0}^{R_1}\ldots\sum_{p_n=0}^{R_n}c_{p_1\ldots p_n}\prod_{i=1}^n\left[\int_{-\infty}^{\infty}l_{p_i}(\tau_i)f(t-\tau_i)d\tau_i\right]
\end{equation}
By denoting
\begin{equation}\label{x1}
  x_i(t)=\int_{-\infty}^{\infty}l_{p_i}(\tau_i)f(t-\tau_i)d\tau_i
\end{equation}
Eq.~\ref{y3} becomes
\begin{equation}\label{y4}
  y_n(t)=\sum_{p_1=0}^{R_1}\ldots\sum_{p_n=0}^{R_n}c_{p_1\ldots p_n}\prod_{i=1}^{n}x_i(t)
\end{equation}

\begin{figure}[H]
  \centering
 \includegraphics[width=5 in]{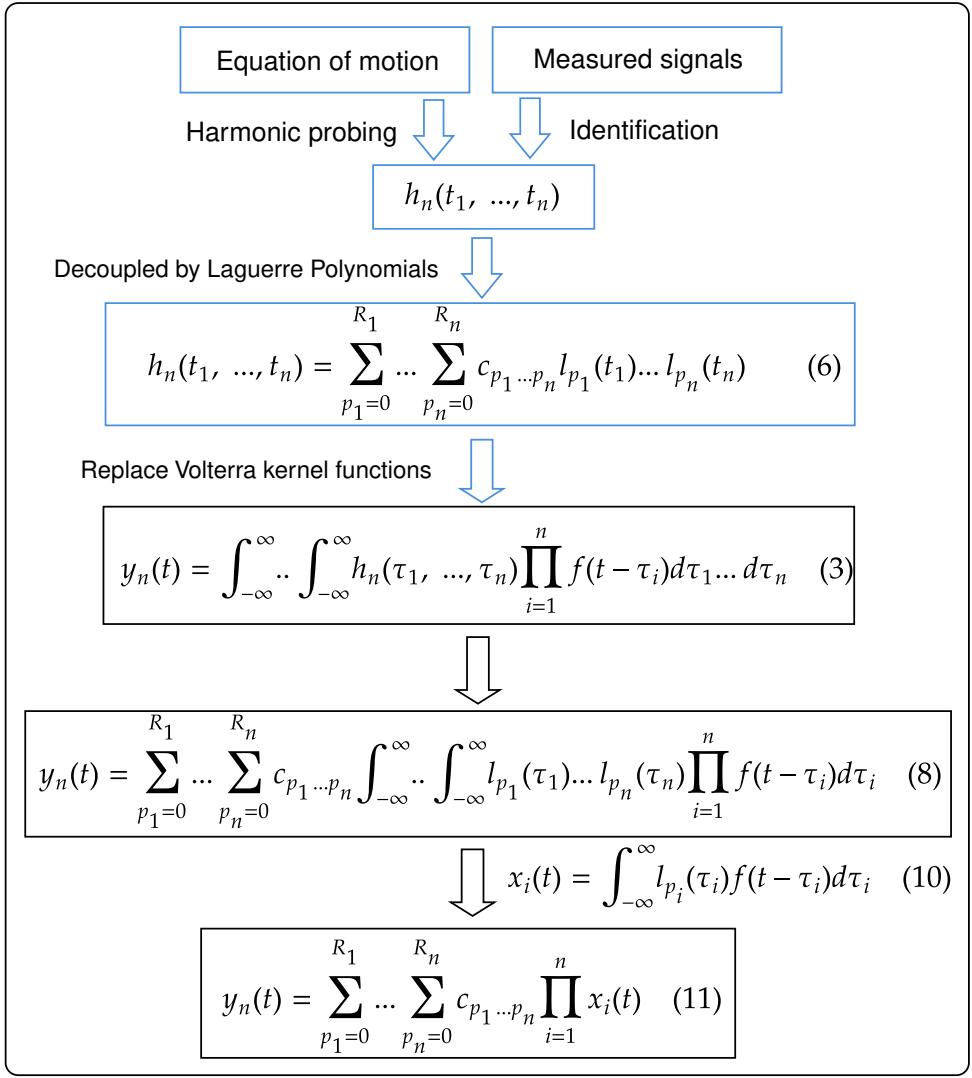}\\
\caption{Procedure to compute the response by a combination of Volterra series and Laguerre polynomials}\label{diag1_volterra}
\end{figure}

The above procedure to compute the nonlinear response by a combination of Volterra series and Laguerre polynomials is schematically shown in Fig.~\ref{diag1_volterra}. Volterra kernel functions $h_{n}(t_{1},\ldots,t_{n})$ can be obtained by either an equation of motion or measured input--output signals.
To derive a closed-form solution of the response, we must obtain a closed-form solution of $x_i(t)$ first. In the following presentation, a closed-form solution of the aforementioned $x_i(t)$ and $y_n(t)$ is derived by using the pole-residue method.

\section{Pole-residue method for calculating $x_i(t)$ and $y_n(t)$}
Performing the Laplace transform of $x_i(t)$ in Eq.~\ref{x1} yields
\begin{equation}\label{x1_s}
  \tilde{x}_i(s)=\tilde{l}_{p_i}(s)\tilde{f}(s)
\end{equation}
where
\begin{equation}\label{laguerre_s}
  \tilde{l}_{p_i}(s)=\mathcal{L}[l_{p_i}(t)]=\sum_{k=0}^{p_i}\frac{b_{p_i}(k)}{(s+a_i)^{k+1}}
\end{equation}
in which
\begin{equation}\label{coef}
 b_{p_i}(k) =\frac{(-1)^{p_i-k}p_i!(2a_i)^{k+0.5}}{k!(p_i-k)!}
\end{equation}
Eq.~\ref{laguerre_s} includes a single pole and several higher-order poles. For $k=0$, $-a_i$ is a single pole, and $b_{p_i}(0)$ is a corresponding coefficient, namely, the residue. For $k>0$, $-a_i$ are higher-order poles, and $b_{p_i}(k)$ are corresponding coefficients.

For an irregular excitation signal $f(t)$ of a finite duration of $T$, it can always be approximated into a pole-residue form by using the complex exponential signal decomposition method---Prony-SS~\cite{hu2013signal}:
\begin{equation}\label{force}
   f(t)=\sum _{\ell=1}^{N_\ell} \alpha_\ell\,  \exp(\lambda_\ell t)    ~~~~~~ 0 \leq t < T
\end{equation}
where $N_\ell$ is the number of components; $\alpha_{\ell}$ and $\lambda_{\ell}$ are constant coefficients, which either are real numbers or occur in complex conjugate pairs. We define $\lambda_{\ell}=-\delta_{\ell}+\mathrm{i}\Omega_{\ell}$, where $\Omega_{\ell}$ is the excitation frequency and $\delta_{\ell}$ is the damping factor of the $\ell^{th}$ component. We denote $\alpha_{\ell}=A_{\ell}e^{\mathrm{i}\theta_{\ell}}$, where $A_{\ell}$ is the amplitude and $\theta_{\ell}$ is the sinusoidal initial phase in radians. Taking the Laplace transform of Eq.~\ref{force} yields
\begin{equation}\label{force_s}
\tilde{f}(s)= \sum_{\ell=1}^{N_\ell}  \frac{\alpha_\ell}{s-\lambda_\ell}
\end{equation}
Note that the concept of the Prony-SS method is similar to that of a principal component method. A smooth excitation usually requires just several terms to achieve a good approximation. For high irregular loadings, including more terms would achieve a better approximation.

Substituting Eqs.~\ref{laguerre_s} and \ref{force_s} into Eq.~\ref{x1_s} yields
\begin{equation}\label{x1_s2}
  \tilde{x}_i(s)=\left[\sum_{k=0}^{p_i}\frac{b_{p_i}(k)}{(s+a_i)^{k+1}}\right]\left(\sum_{\ell=1}^{N_\ell}  \frac{\alpha_\ell}{s-\lambda_\ell}\right)
\end{equation}
Expressing $\tilde{x}_i(s)$ in its pole-residue form~\cite{kreyszig2010advanced,craig2006fundamentals} yields
\begin{equation}\label{x1_s3}
  \tilde{x}_i(s)=\sum_{k=0}^{p_i}\frac{\beta_{p_i,k}}{(s+a_i)^{k+1}}+\sum_{\ell=1}^{N_\ell}\frac{\gamma_{p_i,\ell}}{s-\lambda_\ell}
\end{equation}
where $\lambda_\ell$ are simple poles, and the corresponding residues are easily obtained by
\begin{equation}\label{eq:residue_last}
  \gamma_{p_i,\ell}=\lim\limits_{s \rightarrow \lambda_\ell} (s-\lambda_\ell)  \tilde{x}_i(s)=\sum_{k=0}^{p_i}\frac{\alpha_\ell b_{p_i}(k)}{(\lambda_\ell+a_i)^{k+1}}
\end{equation}
and $-a_i$ are higher--order poles, and the corresponding coefficients are firstly derived as:
\begin{equation}\label{eq:residue_first}
 \beta_{p_i,k}=\lim\limits_{s \rightarrow -a}\frac{1}{(p_i-k)!}\frac{d^{p_i-k}}{ds^{p_i-k}}\left[(s+a_i)^{p_i+1}\tilde{x}_i(s)\right]=\sum_{\ell=1}^{N_\ell}\sum_{q=1}^{p_i-k+1}\frac{(-1)^{q-1}\alpha_{\ell}b_{p_i}(k+q-1)}{(-a_i-\lambda_\ell)^q}
\end{equation}
By taking the inverse Laplace transform of Eq.~\ref{x1_s3}, a closed-form solution is obtained:
\begin{equation}\label{x3}
 x_i(t)=\sum_{k=0}^{p_i}\frac{\beta_{p_i,k}}{k!}t^{k}e^{-a_it}+\sum_{\ell=1}^{N_\ell}\gamma_{p_i,\ell}e^{\lambda_\ell t}
\end{equation}
Substituting Eqs.~\ref{y4} and ~\ref{x3} into Eq.~\ref{y1} yields
\begin{equation}\label{y4_final}
  y(t)=\sum_{n=1}^N\sum_{p_1=0}^{R_1}\ldots\sum_{p_n=0}^{R_n}c_{p_1\ldots p_n}\prod_{i=1}^{n}\left(\sum_{k=0}^{p_i}\frac{\beta_{p_i,k}}{k!}t^{k}e^{-a_it}+\sum_{\ell=1}^{N_\ell}\gamma_{p_i,\ell}e^{\lambda_\ell t}\right)
\end{equation}

Theoretically speaking, the proposed method for deriving the closed-form solution of the nonlinear response is applicable to any order of the Volterra series. For practical engineering, usually only the first several order responses dominate. By setting up $N=2$, Eq.~\ref{y4_final} can be simplified into three components:
\begin{equation}\label{y5}
\begin{aligned}
  y(t)=y_s(t)+y_c(t)+y_f(t)
\end{aligned}
\end{equation}
where the natural response, which is only related to system poles, is given by
\begin{equation}\label{y51_natural}
\begin{aligned}
  y_s(t)=\sum_{p_1=0}^{R_1}\sum_{k=0}^{p_1}\frac{c_{p_1}\beta_{p_1,k}}{k!}t^{k}e^{-a_1t}+\sum_{p_1=0}^{R_1}\sum_{p_2=0}^{R_2}\sum_{k=0}^{p_1}\sum_{m=0}^{p_2}\frac{c_{p_1 p_2}\beta_{p_1,k}\beta_{p_2,m}}{k!m!}t^{(k+m)}e^{-(a_1+a_2)t}
\end{aligned}
\end{equation}
and the cross response, which is related to both system poles and excitation poles, is given by
\begin{equation}\label{y51_cross}
\begin{aligned}
  y_c(t)&=\sum_{p_1=0}^{R_1}\sum_{p_2=0}^{R_2}\sum_{k=0}^{p_1}\sum_{j=1}^{N_\ell}\frac{c_{p_1 p_2}\beta_{p_1,k}\gamma_{p_2,j}t^{k}}{k!}e^{(-a_1+\lambda_j)t}\\
  &+\sum_{p_1=0}^{R_1}\sum_{p_2=0}^{R_2}\sum_{m=0}^{p_2}\sum_{\ell=1}^{N_\ell}\frac{c_{p_1 p_2}\beta_{p_2,m}\gamma_{p_1,\ell}t^{m}}{m!}e^{(-a_2+\lambda_\ell)t}
\end{aligned}
\end{equation}
and the forced response, which is related only to excitation poles, is given by
\begin{equation}\label{y51_forced2}
\begin{aligned}
  y_f(t)=\sum_{p_1=0}^{R_1}\sum_{\ell=1}^{N_\ell}c_{p_1}\gamma_{p_1,\ell}e^{\lambda_\ell t}+\sum_{p_1=0}^{R_1}\sum_{p_2=0}^{R_2}\sum_{\ell=1}^{N_\ell}\sum_{j=1}^{N_\ell}c_{p_1 p_2}\gamma_{p_1,\ell}\gamma_{p_2,j}e^{(\lambda_\ell+\lambda_j) t}
\end{aligned}
\end{equation}
The first term in Eq.~\ref{y51_forced2} is the first-order forced response governed by the excitation frequency, \textit{i.e.}, the imaginary part of the pole $\lambda_\ell$. The second term corresponds to the second-order nonlinear forced response, which includes the sum frequency and difference frequency responses governed by $\lambda_\ell+\lambda_j$. Eq.~\ref{y51_forced2} straightforwardly offers visible information about the possible nonlinear vibrations by the cooperation of excitation frequencies.

Particularly, consider a sinusoidal excitation $f(t)=\sin \omega_r t$, which can be expressed as $f(t)=\gamma e^{\lambda t}+\gamma^* e^{\lambda^*t}$, where $\gamma=-0.5\mathrm{i}$ and $\lambda=\mathrm{i}\omega_r$. Substituting these values into Eq.~\ref{y51_forced2}, the second term of Eq.~\ref{y51_forced2} is simplified as
\begin{equation}\label{y51_forced3}
\begin{aligned}
  y_{f_2}(t)=\sum_{p_1=0}^{R_1}\sum_{p_2=0}^{R_2}0.5c_{p_1 p_2}-\sum_{p_1=0}^{R_1}\sum_{p_2=0}^{R_2}0.5c_{p_1 p_2}\cos 2\omega_rt
\end{aligned}
\end{equation}
where the first term is the difference frequency response, and the second term is the sum frequency response.

\section{Numerical studies}
In practical engineering, some systems have an accurate equation of motion. Additionally, some systems have difficulty constructing their equations of motion because of complex nonlinear dynamic behaviours and uncertain system parameters. In this article, a system with a known equation of motion is called a known system, and a system with an unknown equation of motion is called an unknown system for simplicity. In this section, two numerical studies are presented. The first study verifies the proposed method using a known nonlinear oscillator, and the second study demonstrates the applicability of the proposed method to an unknown system.
Throughout the numerical studies, the unit system is the metre--kilogramme--second (MKS) system; for conciseness, explicit units for quantities are omitted.

\subsection{A known nonlinear system}\label{example1}
This study chooses a nonlinear oscillator written as:
\begin{equation}\label{eq:nonlinear_SDOF}
  m\ddot{y}+c\dot{y}+k_1y+k_2y^{2}+k_3y^{3}=f(t)
\end{equation}
where mass $m=1$, damping $c=1$, linear stiffness $k_1=10$, quadratic stiffness $k_2=20$ and cubic stiffness $k_3=20$. It is a case that has been studied in a previously published article~\cite{kim2014time}. The linear natural frequency of the system $\omega_0=\sqrt{k_1/m}=3.16$ and the damping ratio $\zeta=c/(2m\omega_0)=15.8\%$. This kind of oscillator occurs in many engineering problems, such as a model of fluid resonance in a narrow gap between large vessels~\cite{song2021}. In the model, $k_1y$ represents the linear restoring force of the fluid, and $k_2y^{2}$ and $k_3y^{3}$ are respectively the quadratic and cubic nonlinear restoring forces of the fluid.

\subsubsection{Volterra kernel functions}
\begin{figure}[H]
  \centering
 \includegraphics[width=4 in]{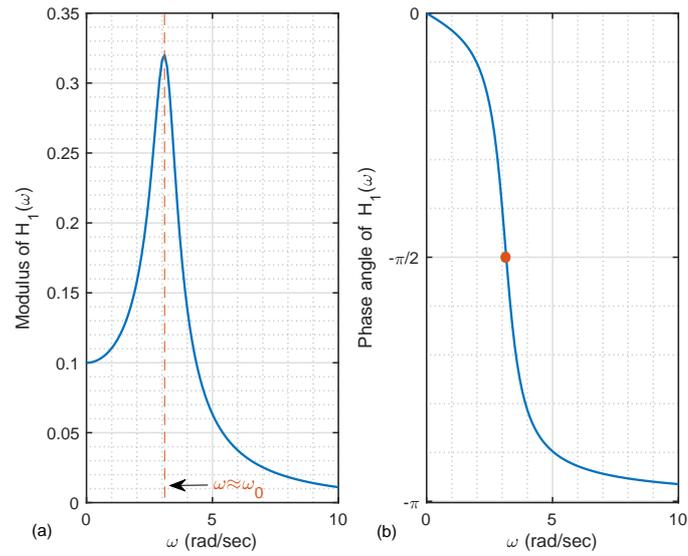}\\
\caption{Linear frequency response function: (a) Modulus of $H_1(\omega)$, (b) phase angle of $H_1(\omega)$}\label{LFRF}
\end{figure}

\begin{figure}[H]
  \centering
 \includegraphics[width=4 in]{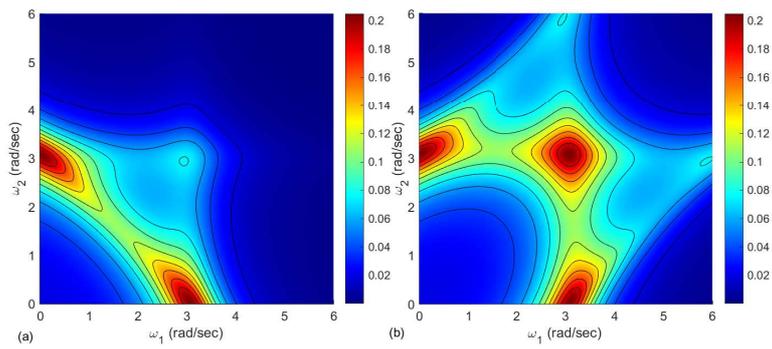}\\
\caption{Modulus of quadratic frequency response function: (a) sum frequency, (b) difference frequency}\label{QFRF}
\end{figure}

\begin{figure}[H]
  \centering
 \includegraphics[width=4 in]{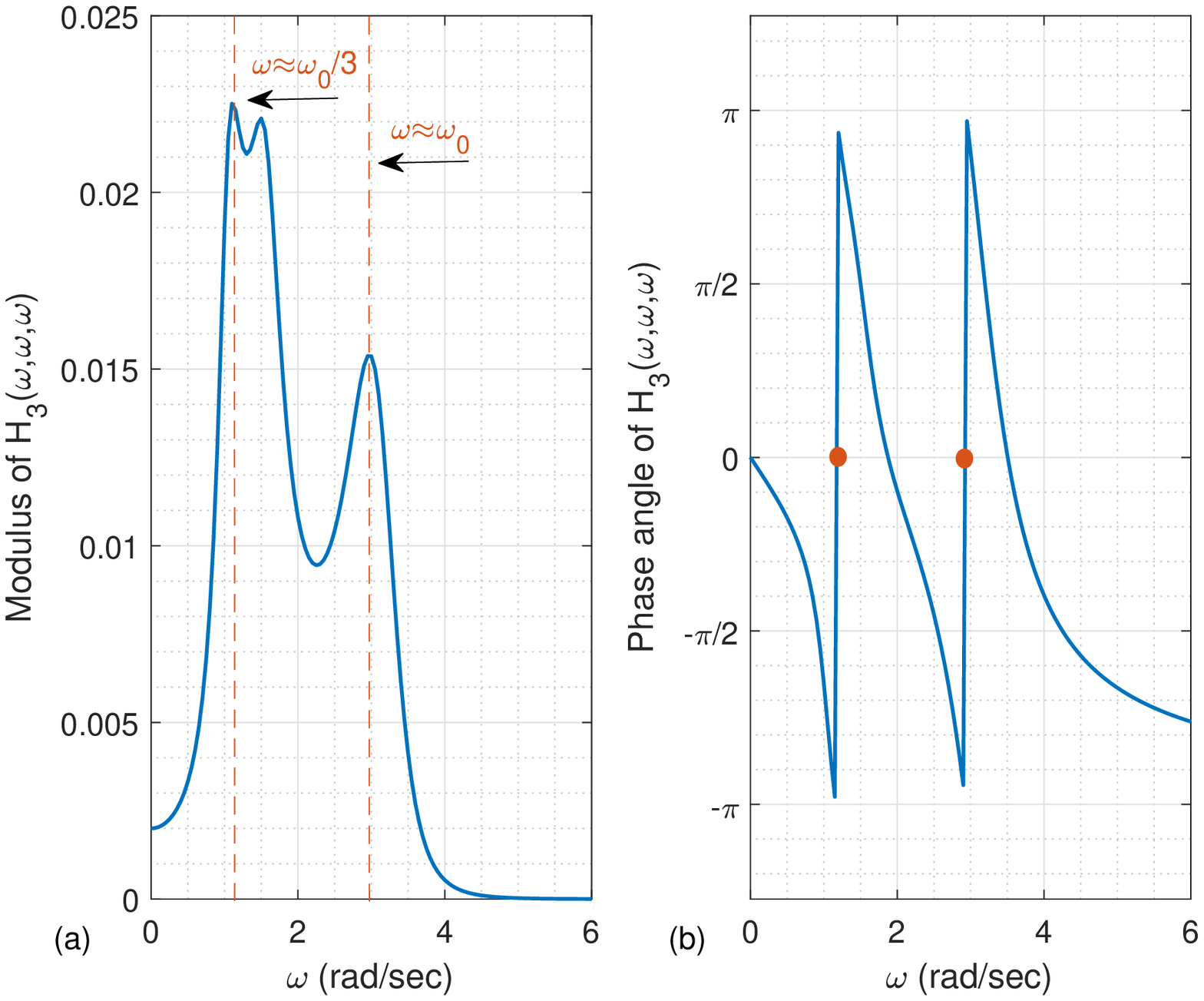}\\
\caption{Diagonal term of the sum frequency cubic frequency response function: (a) modulus of $H_3(\omega,\omega,\omega)$, (b) phase angle of $H_3(\omega,\omega,\omega)$}\label{CFRF_sum}
\end{figure}

\begin{figure}[H]
  \centering
 \includegraphics[width=4 in]{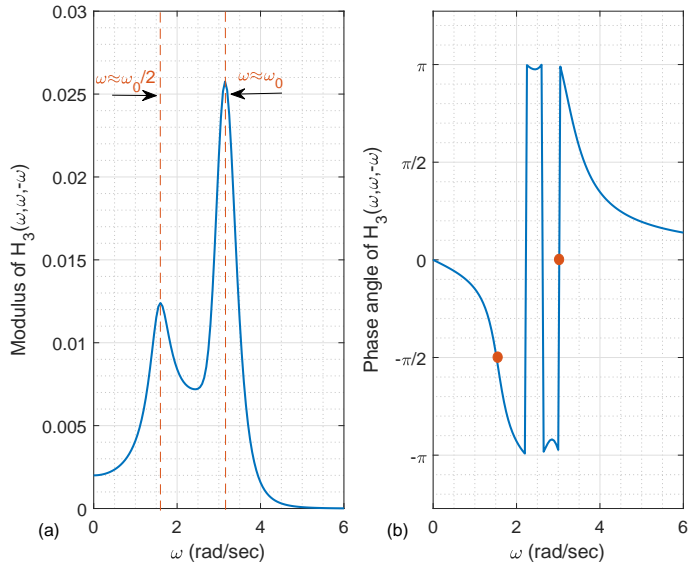}\\
\caption{Diagonal term of the difference frequency cubic frequency response function: (a) modulus of $H_3(\omega,\omega,-\omega)$, (b) phase angle of $H_3(\omega,\omega,-\omega)$}\label{CFRF_diff}
\end{figure}

Generally, the first several order responses dominate the total response of a system. Hence, the order of the Volterra series in Eq.~\ref{y4_final} is chosen to be 3, namely, $N=3$. For computing the first three order responses from Eq.~\ref{y4_final}, the first three order Volterra kernel functions need to be known. Since Volterra kernel functions and corresponding frequency response functions are related by a specific Fourier transform pair, we can first write the first three orders of frequency response functions directly from Eq.~\ref{eq:nonlinear_SDOF}. Then, Volterra kernel functions are obtained by the inverse Fourier transform. Based on the harmonic probing algorithm~\cite{bedrosian1971output,chatterjee2010parameter}, the linear frequency response function (LFRF) $H_1(\omega)$, the quadratic frequency response function (QFRF) $H_2(\omega_1,\omega_2)$ and the cubic frequency response function (CFRF) $H_3(\omega_1,\omega_2,\omega_3)$ are analytically given by:
\begin{equation}\label{eq:LFRF}
  H_1(\omega)=\frac{1}{-m\omega^2+\mathrm{i}c\omega+k_1}, ~-\infty<\omega<\infty
\end{equation}
\begin{equation}\label{eq:QFRF}
  H_2(\omega_1,\omega_2)=-k_2H_1(\omega_1)H_1(\omega_2)H_1(\omega_1+\omega_2), ~-\infty<\omega_1,\omega_2<\infty
\end{equation}
and
\begin{equation}\label{eq:CFRF}
\begin{aligned}
  &H_3(\omega_1,\omega_2,\omega_3)=-\left\{\frac{k_2}{3}\left[H_1(\omega_1)H_2(\omega_2,\omega_3)+H_1(\omega_2)H_2(\omega_1,\omega_3)+H_1(\omega_3)H_2(\omega_1,\omega_2)\right]\right.\\
  &\left.+k_3H_1(\omega_1)H_1(\omega_2)H_1(\omega_3)\right\}H_1(\omega_1+\omega_2+\omega_3), ~-\infty<\omega_1, \omega_2, \omega_3<\infty
\end{aligned}
\end{equation}
Figures~\ref{LFRF}-\ref{CFRF_diff} show $H_1(\omega)$, $H_2(\omega_1,\omega_2)$ and $H_3(\omega_1,\omega_2,\omega_3)$, respectively, which agree well with those reported in Ref.~\cite{kim2014time}. As expected, the modulus of $H_1(\omega)$ in Fig.~\ref{LFRF} peaks near the linear natural frequency $\omega_0$, and the phase angle decreases monotonically from 0 to -$\pi$ with increasing frequency. Figure~\ref{QFRF} shows the sum frequency QFRF, where the energy converges along the line of $\omega_1+\omega_2\approx\omega_0$. Therefore, when the sum frequency of a two-tone excitation equals the linear resonant frequency, the second-order response may reach its maximum. Additionally, those pairs of excitations in line $\omega_1+\omega_2\approx\omega_0$ may produce non-negligible vibration magnitudes due to second-order nonlinear effects. For the difference frequency QFRF in Fig.~\ref{QFRF}(b), the energy converges along two main lines, \textit{i.e.}, $\omega_1\approx\omega_0$ and $\omega_2\approx\omega_0$. Figures~\ref{CFRF_sum} and \ref{CFRF_diff} show moduli of $H_3(\omega,\omega,\omega)$ and $H_3(\omega,\omega,-\omega)$, which are diagonal terms of the sum frequency CFRF and the difference frequency CFRF, respectively. While the modulus of $H_3(\omega,\omega,\omega)$ peaks near $\omega\approx\omega_0/3$ and $\omega_0$, that of $H_3(\omega,\omega,-\omega)$ peaks near $\omega\approx\omega_0$ with a small hump around $\omega\approx\omega_0/2$. Values at $\omega\approx\omega_0/3$ and $\omega_0/2$ may be magnified by higher-order stiffness terms in Eq.~\ref{eq:nonlinear_SDOF}.

\begin{figure}[H]
 \centering
\includegraphics[width=4in]{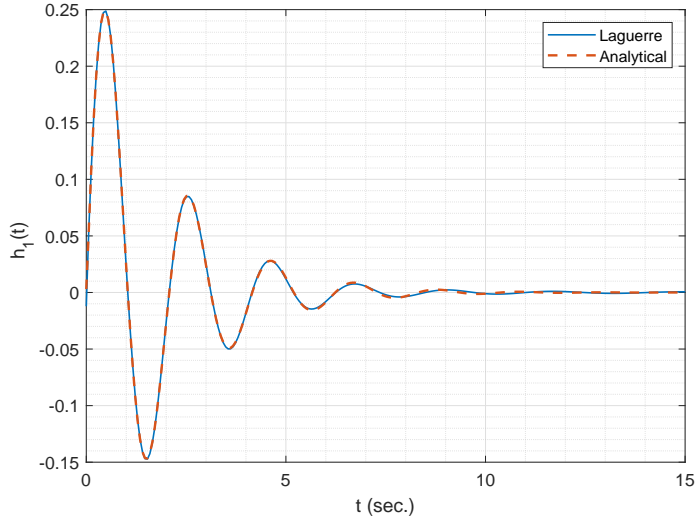}
\caption{Comparison of $h_1(t)$ based on the analytical and reconstructed by Laguerre polynomials}\label{LIRF}
\end{figure}

\begin{figure}[H]
\centering
\begin{minipage}[t]{0.5\linewidth}
\centering
\includegraphics[width=3in]{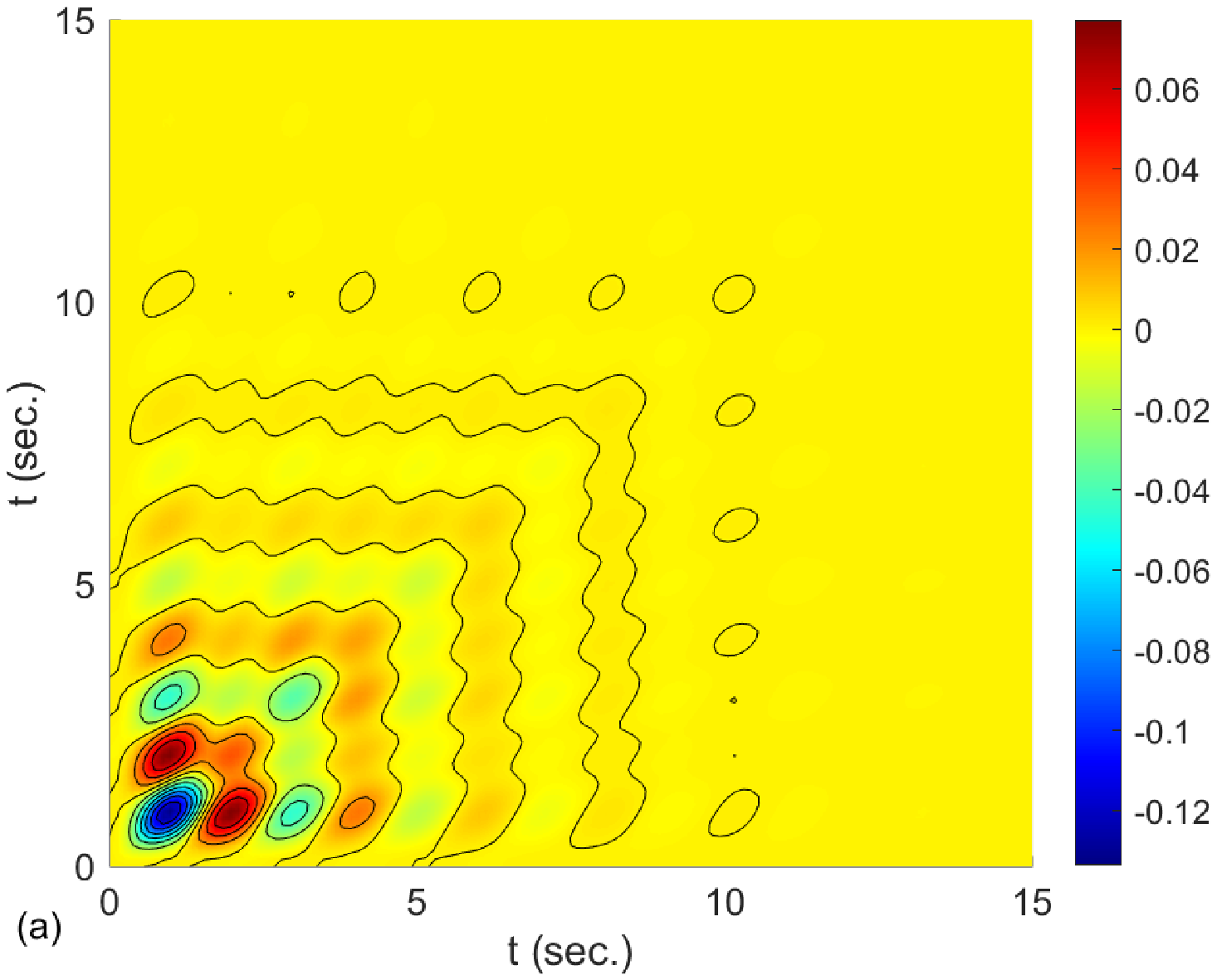}
\end{minipage}
\begin{minipage}[t]{0.5\linewidth}
\centering
\includegraphics[width=3in]{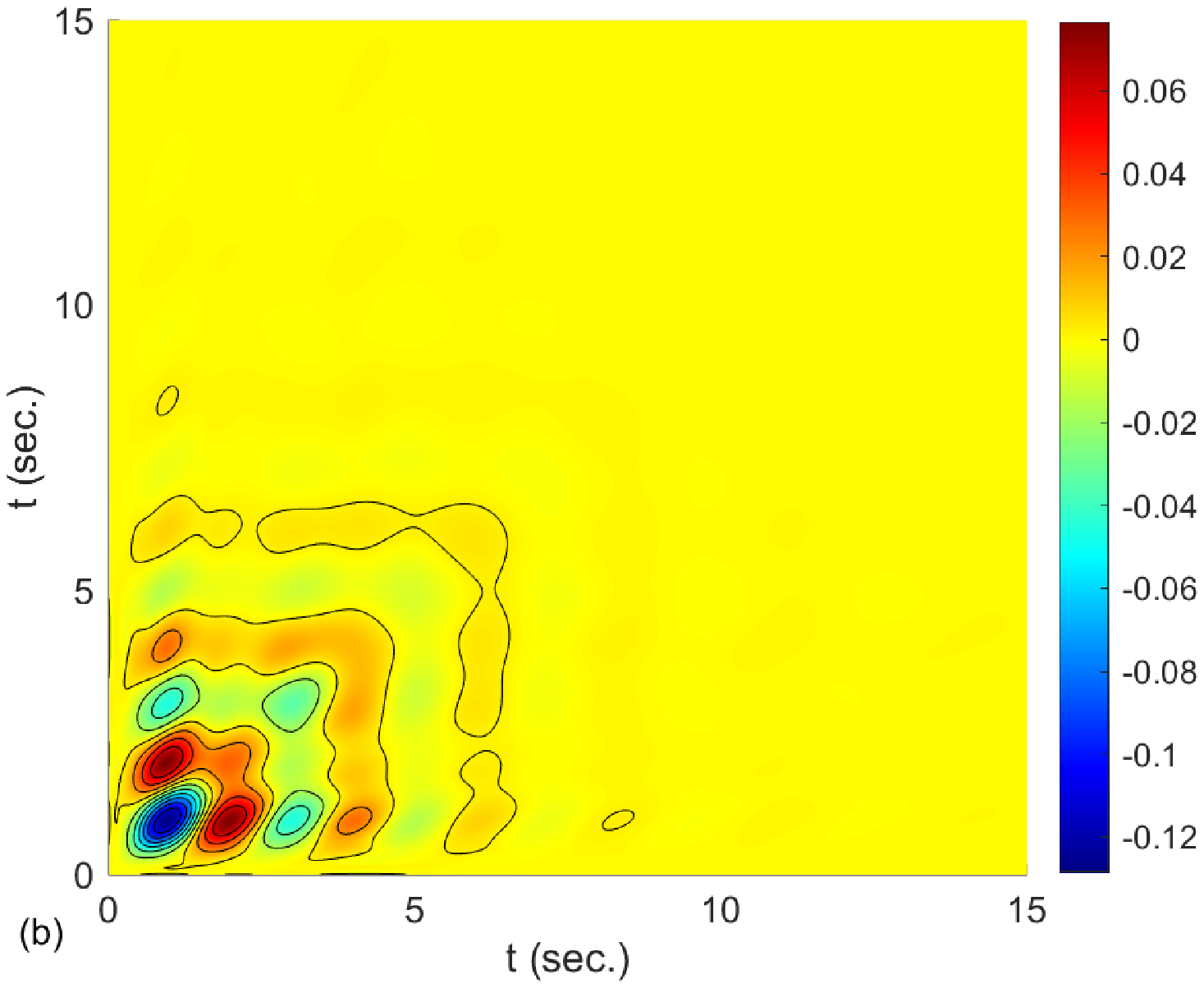}
\end{minipage}
\centering
\caption{Second-order Volterra kernel functions: $h_2(t_1,t_2)$: (a) analytical, (b) reconstructed by Laguerre polynomials}\label{QIRF}
\end{figure}

\begin{figure}[H]
\centering
\includegraphics[width=4in]{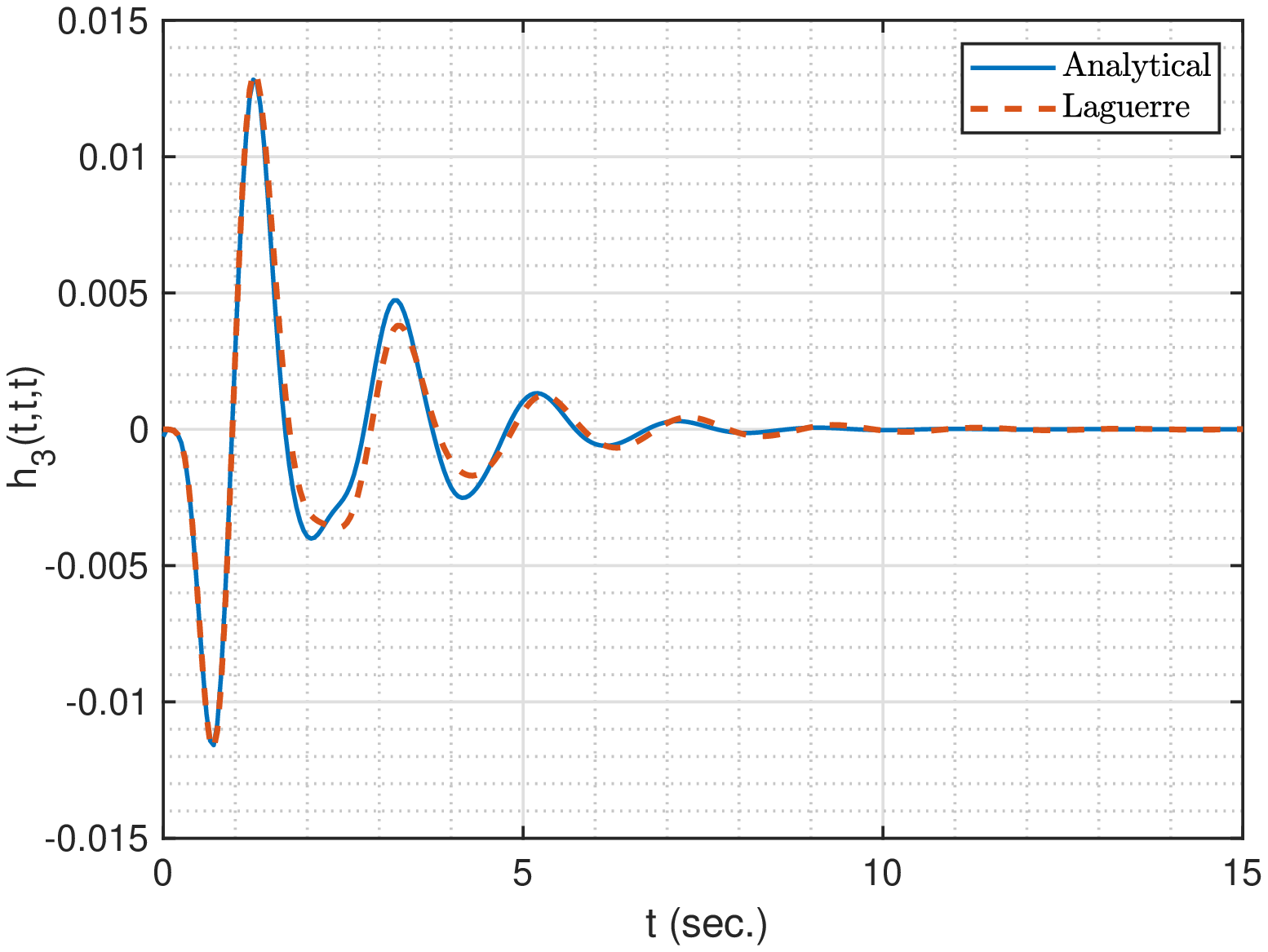}
\caption{Comparison of $h_3(t,t,t)$ based on the analytical and reconstructed by Laguerre polynomials}\label{CIRF_diag}
\end{figure}

\begin{figure}[H]
\centering
\begin{minipage}[t]{0.5\linewidth}
\centering
\includegraphics[width=3in]{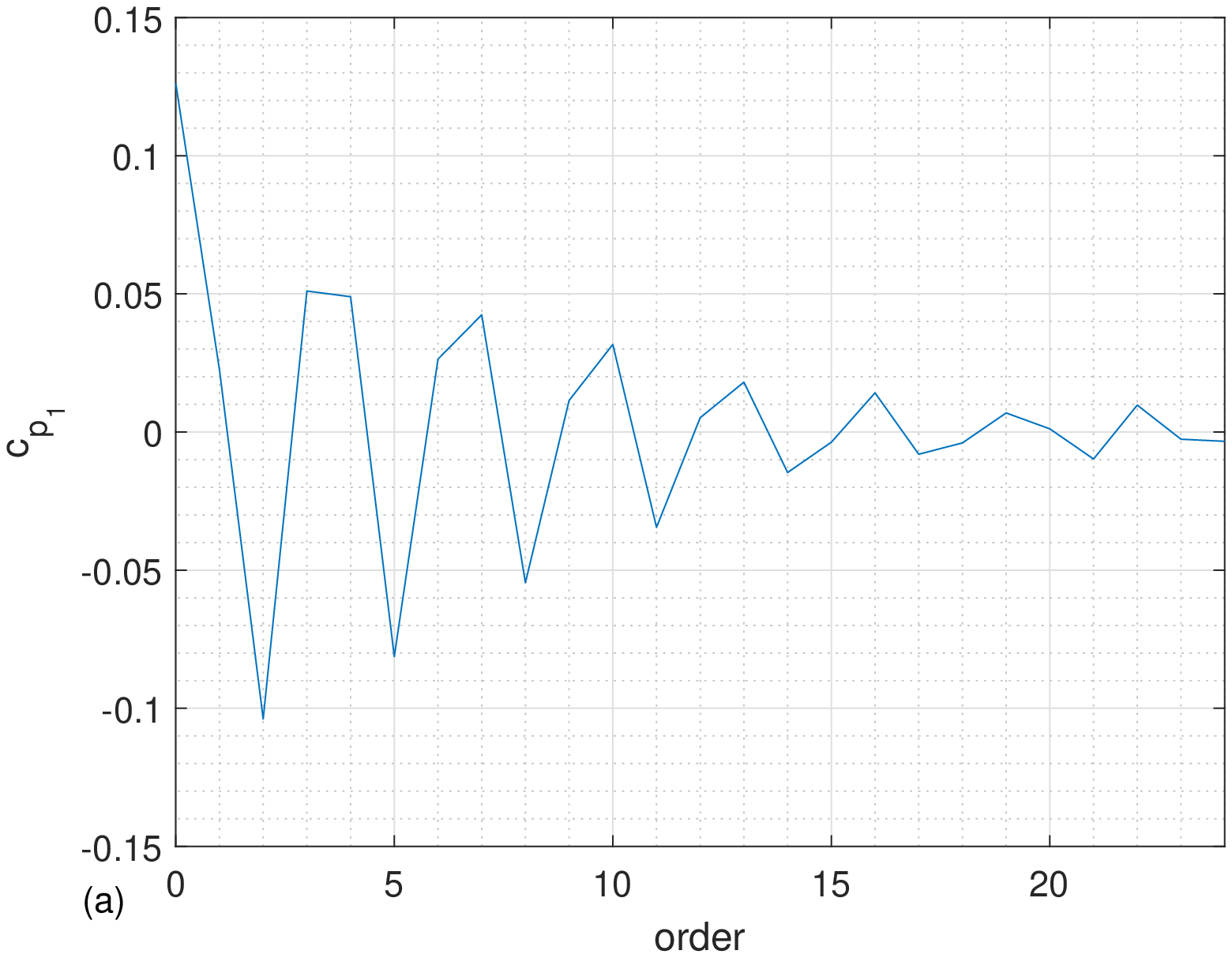}
\end{minipage}
\begin{minipage}[t]{0.5\linewidth}
\centering
\includegraphics[width=3in]{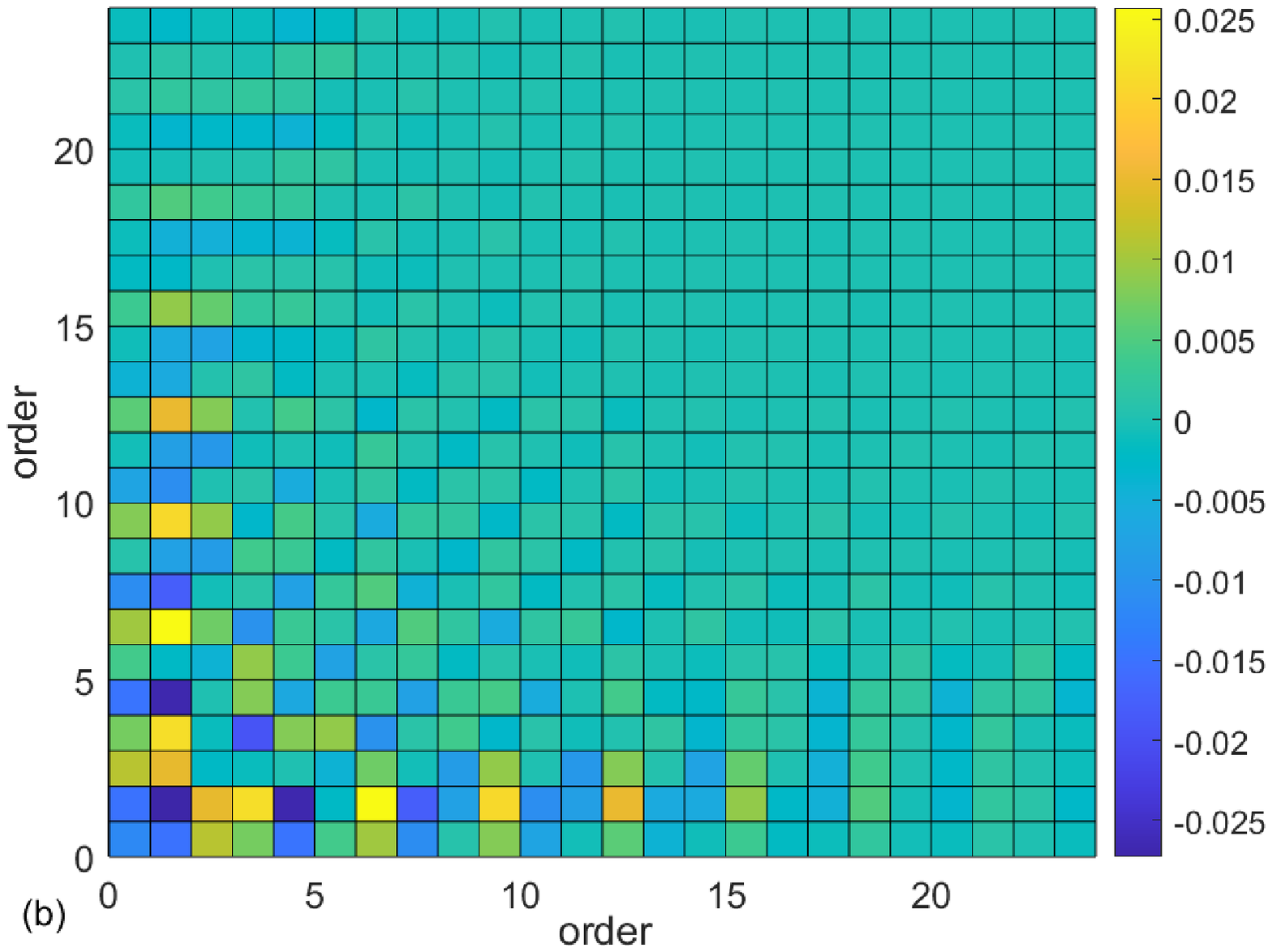}
\end{minipage}
\centering
\caption{Coefficients of first two orders of Volterra kernel functions: (a) $c_{p_1}$, (b) $c_{p_1p_2}$}\label{coef_exp1}
\end{figure}

\begin{figure}[H]
\centering
\includegraphics[width=4in]{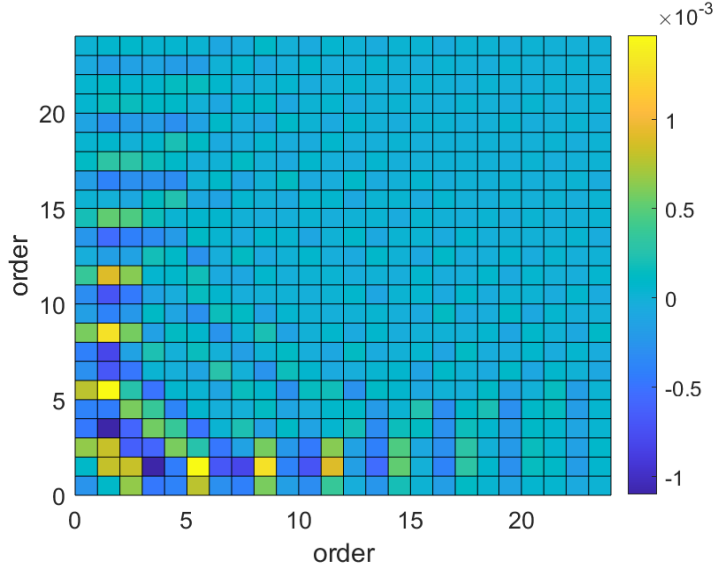}
\caption{Coefficients of the third-order Volterra kernel function: $c_{p_1p_2p_3}$ for $p_3=0$}\label{c31}
\end{figure}

By performing the inverse fast Fourier transform to Eqs.~\ref{eq:LFRF}-\ref{eq:CFRF}, the corresponding linear impulse response function $h_1(t)$, quadratic impulse response function $h_2(t_1,t_2)$ and cubic impulse response function $h_3(t_1,t_2,t_3)$ are obtained. Here, $h_1(t)$ and $h_2(t_1,t_2)$ are plotted in Figs.~\ref{LIRF} and \ref{QIRF}, respectively, and $h_3(t,t,t)$ is shown in Fig.~\ref{CIRF_diag}. In the numerical implementation, Eqs.~\ref{eq:LFRF}-\ref{eq:CFRF} have been utilized with the frequency interval $\Delta\omega=0.1$, number of frequency components $N_n=1025$, and cut-off frequencies $102.4$ and $-102.4$. For decoupling Volterra kernel functions by using Laguerre polynomials, the damping rate and number of Laguerre polynomials for each order Volterra kernel function need to be determined (see Eqs.~\ref{laguerre} and \ref{ht_exp}). In this example, we set $a_1=a_2=a_3=2$ and $R_1=R_2=R_3=24$ because coefficients $c_{p_1\ldots p_n}$ become very small when $R_n>24$, $n=1, 2, 3$. According to Eq.~\ref{eq:coefficient}, the coefficients of the first three order Volterra kernel functions are calculated, which are shown in Figs.~\ref{coef_exp1} and \ref{c31}. For convenience, Fig.~\ref{c31} plots only $c_{p_1p_2p_3}$ for $p_3=0$. With the increase of the order of Laguerre polynomials, coefficients in Figs.~\ref{coef_exp1} and \ref{c31} gradually decrease, which illustrates how the first several orders of Laguerre polynomials dominate all orders of the Volterra kernel function.
With the known Laguerre polynomials and corresponding coefficients, Volterra kernel functions are reconstructed by Eq.~\ref{ht_exp}. For comparison, reconstructed Volterra kernel functions are also plotted in Figs.~\ref{LIRF}-\ref{CIRF_diag}. The reconstructed results agree well with the analytical values, which verifies the accuracy of the decomposition.

\subsubsection{Sinusoidal excitation}\label{sinusoidal}
From Eq.~\ref{eq:nonlinear_SDOF}, we consider
a sinusoidal excitation
\begin{equation}\label{eq:input signal}
     f(t) = A  \sin(\Omega t)
  \end{equation}
where $A$ and $\Omega$ are the amplitude and the frequency, respectively. Five cases of $A$ and $\Omega$ are shown in Table \ref{input parameters}. Excitation frequencies in Cases 1 and 2 are larger than the linear natural frequency ($\omega_0\approx3.16$), those in Case 3 are very close to $\omega_0$, and those in Cases 4 and 5 are smaller than $\omega_0$. All cases have same amplitudes.
The poles of a sinusoidal excitation are $\lambda_{1,2}  =\pm \mathrm{i}\Omega$, and the residues are $\alpha_{1,2} = \mp \mathrm{i}A/2$. Numerical values of excitation poles and residues for different cases are listed in Table \ref{input parameters}.
\begin{table}[h!]
   \begin{center}
   \caption{Parameter values, poles and residues of the sinusoidal excitation}
   \label{input parameters}
   \begin{tabular}{|c|c|c|c|c|}
   \hline
   Case &   $A$ & $\Omega$    &  $\lambda_{1,2}$   & $\alpha_{1,2}$ \\
   \hline
   $1$      & $1$  &  $3\pi$    & $\pm3\pi \mathrm{i}$    & $0.5$\\
   \hline
   $2$      & $1$  &  $2\pi$     & $\pm2\pi \mathrm{i}$    & $0.5$\\
   \hline
   $3$      & $1$  &  $\pi$   & $\pm\pi \mathrm{i}$    & $0.5$\\
   \hline
   $4$      & $1$  &  $0.5\pi$    & $\pm0.5\pi \mathrm{i}$    & $0.5$\\
   \hline
   $5$      & $1$  &  $0.3\pi$    & $\pm0.3\pi \mathrm{i}$    & $0.5$\\
     \hline
   \end{tabular}
   \end{center}
\end{table}

\begin{figure}[H]
\centering
\begin{minipage}[t]{0.5\linewidth}
\centering
\includegraphics[width=3in]{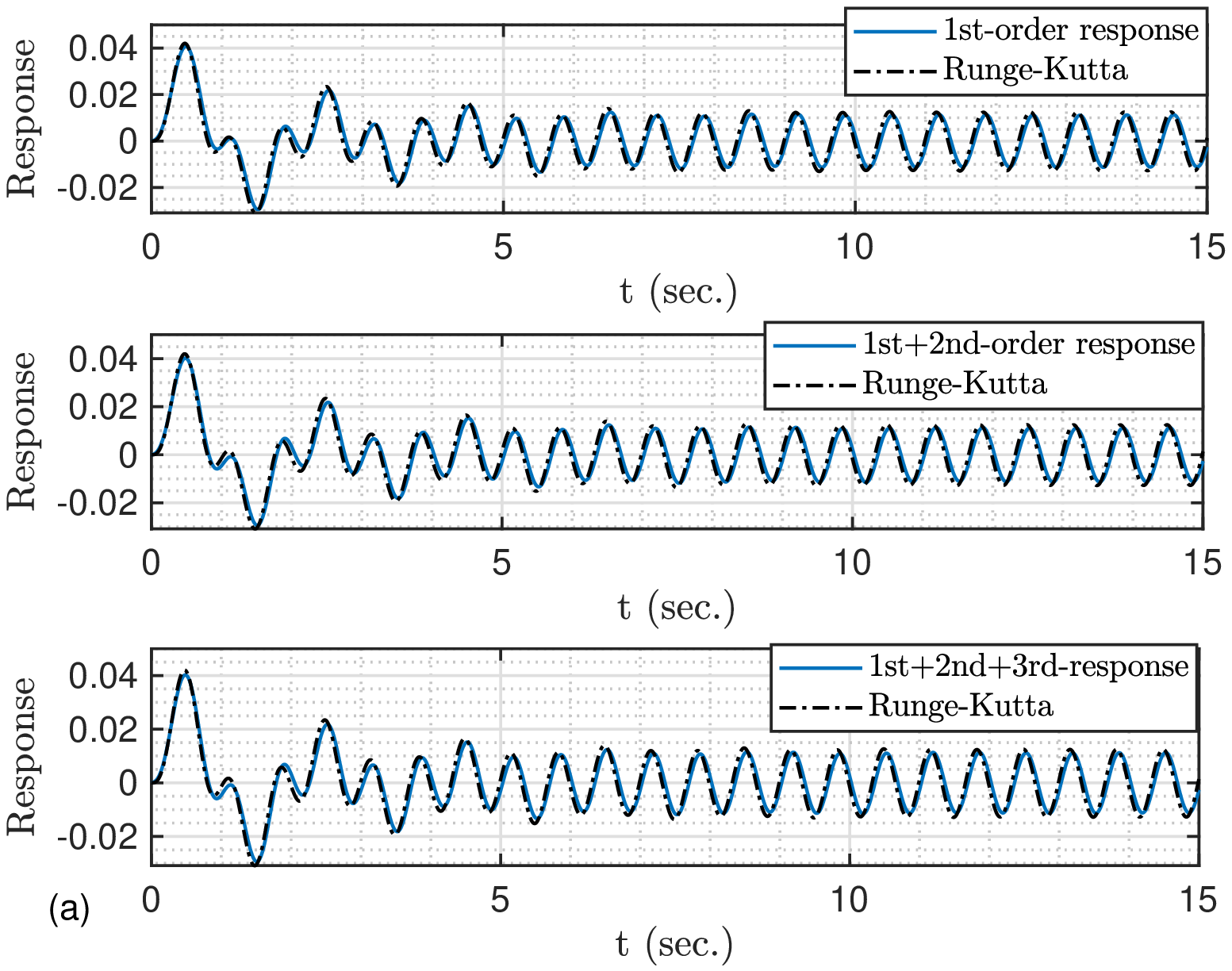}
\end{minipage}
\begin{minipage}[t]{0.5\linewidth}
\centering
\includegraphics[width=3in]{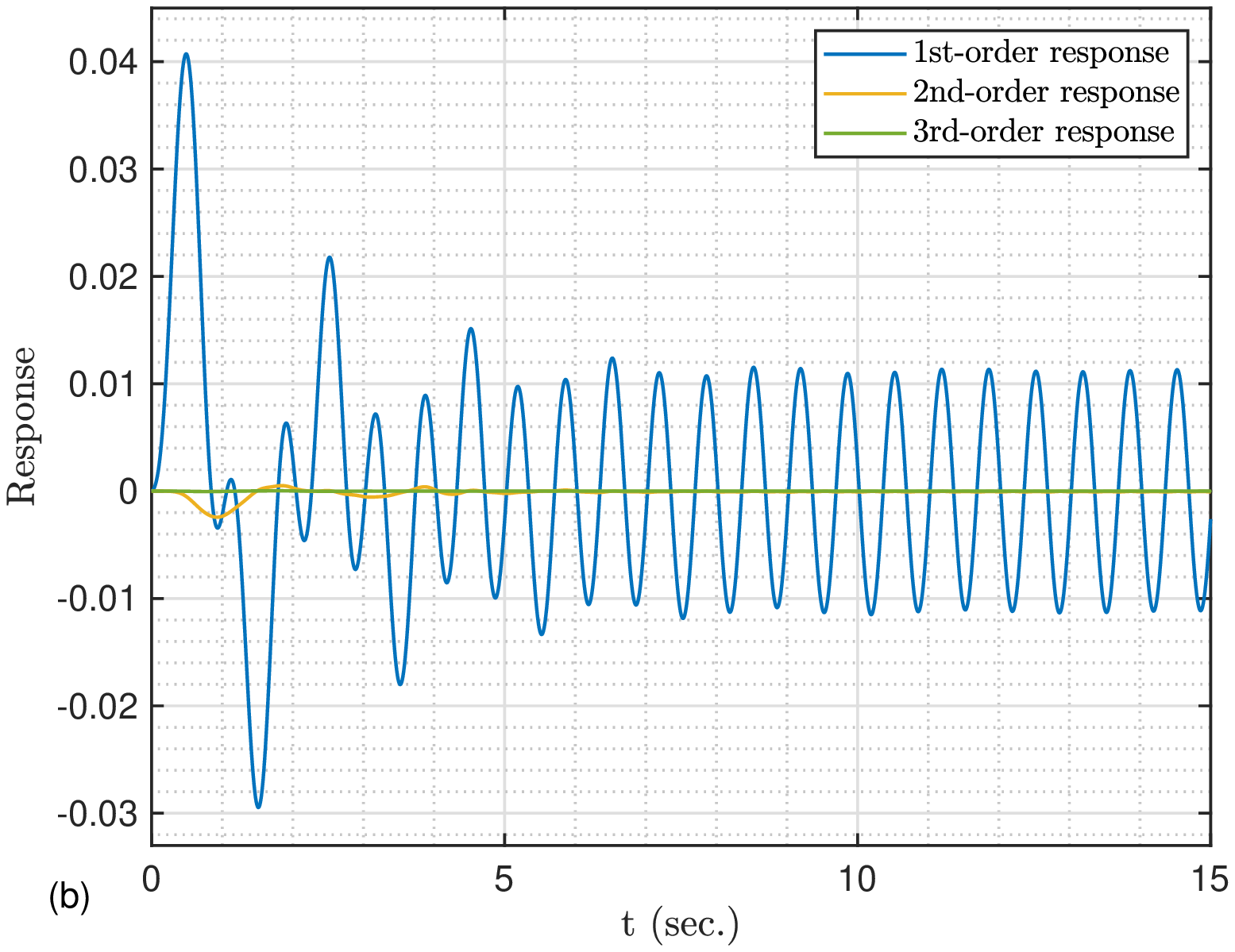}
\end{minipage}
\centering
\caption{Response for Case 1: (a) comparison between the proposed method and Runge--Kutta method, (b) contribution of the three components}\label{exp1_sinu5}
\end{figure}

\begin{figure}[H]
\centering
\begin{minipage}[t]{0.5\linewidth}
\centering
\includegraphics[width=3in]{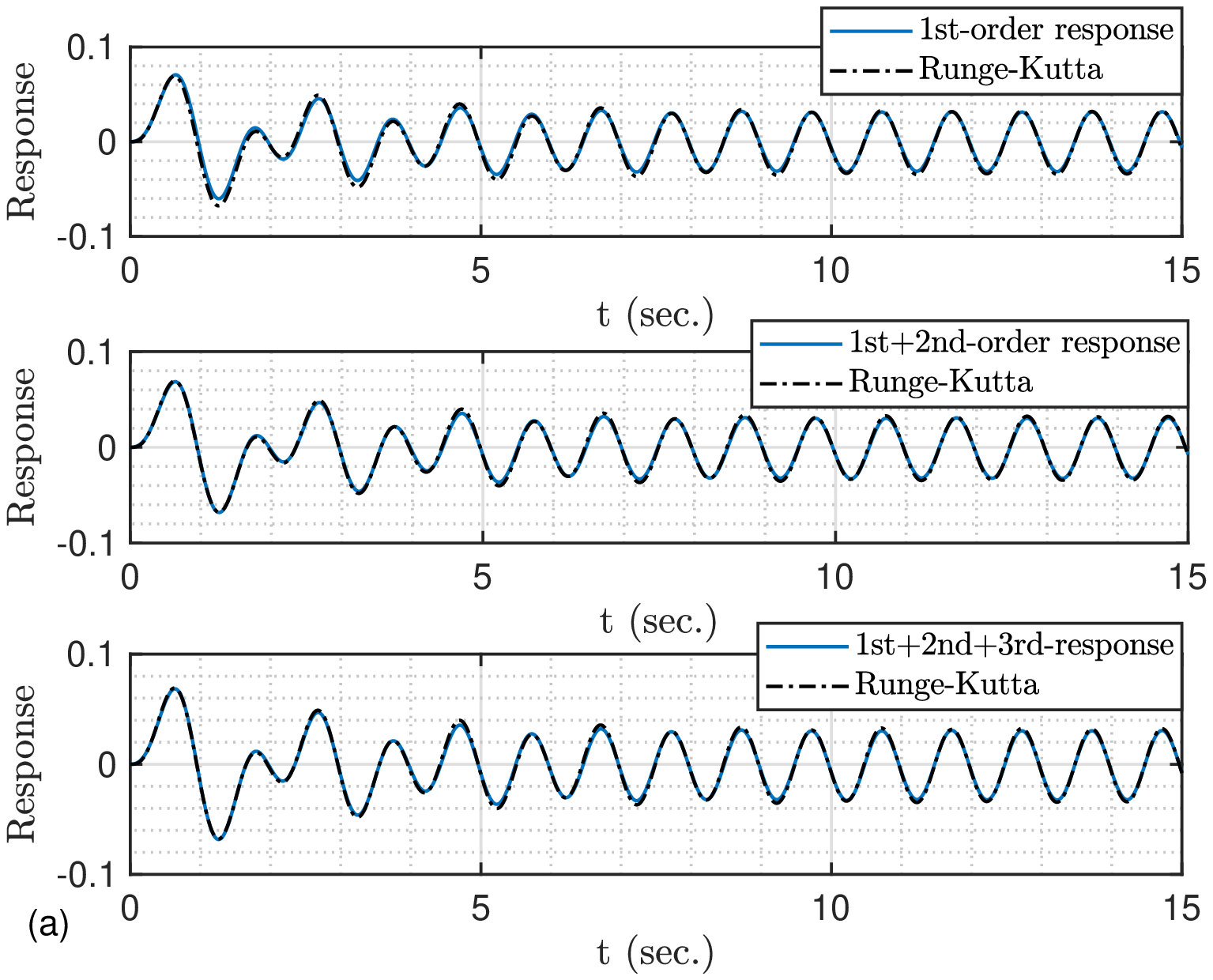}
\end{minipage}
\begin{minipage}[t]{0.5\linewidth}
\centering
\includegraphics[width=3in]{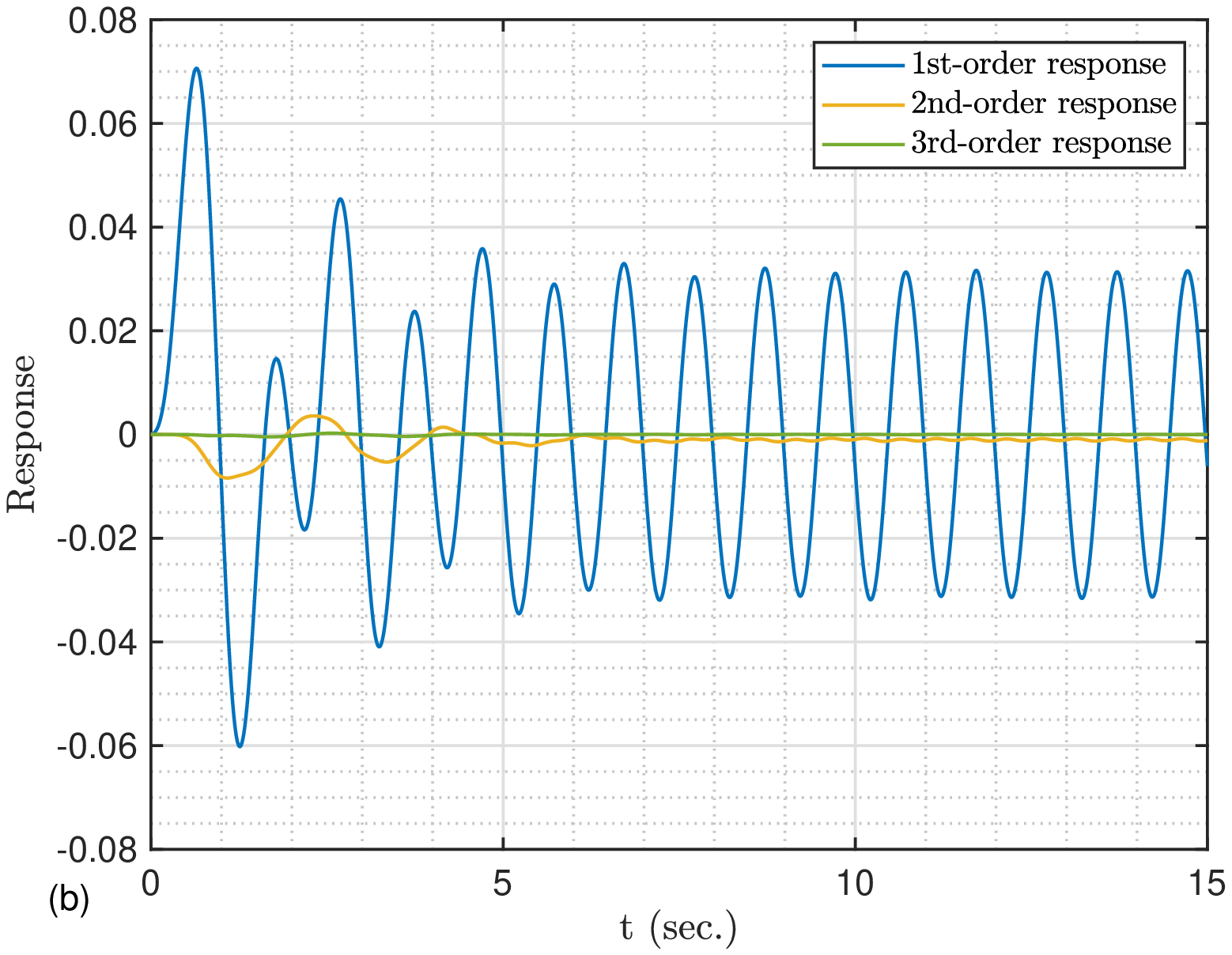}
\end{minipage}
\centering
\caption{Response for Case 2: (a) comparison between the proposed method and Runge--Kutta method, (b) contribution of the three components}\label{exp1_sinu3}
\end{figure}

\begin{figure}[H]
\centering
\begin{minipage}[t]{0.5\linewidth}
\centering
\includegraphics[width=3in]{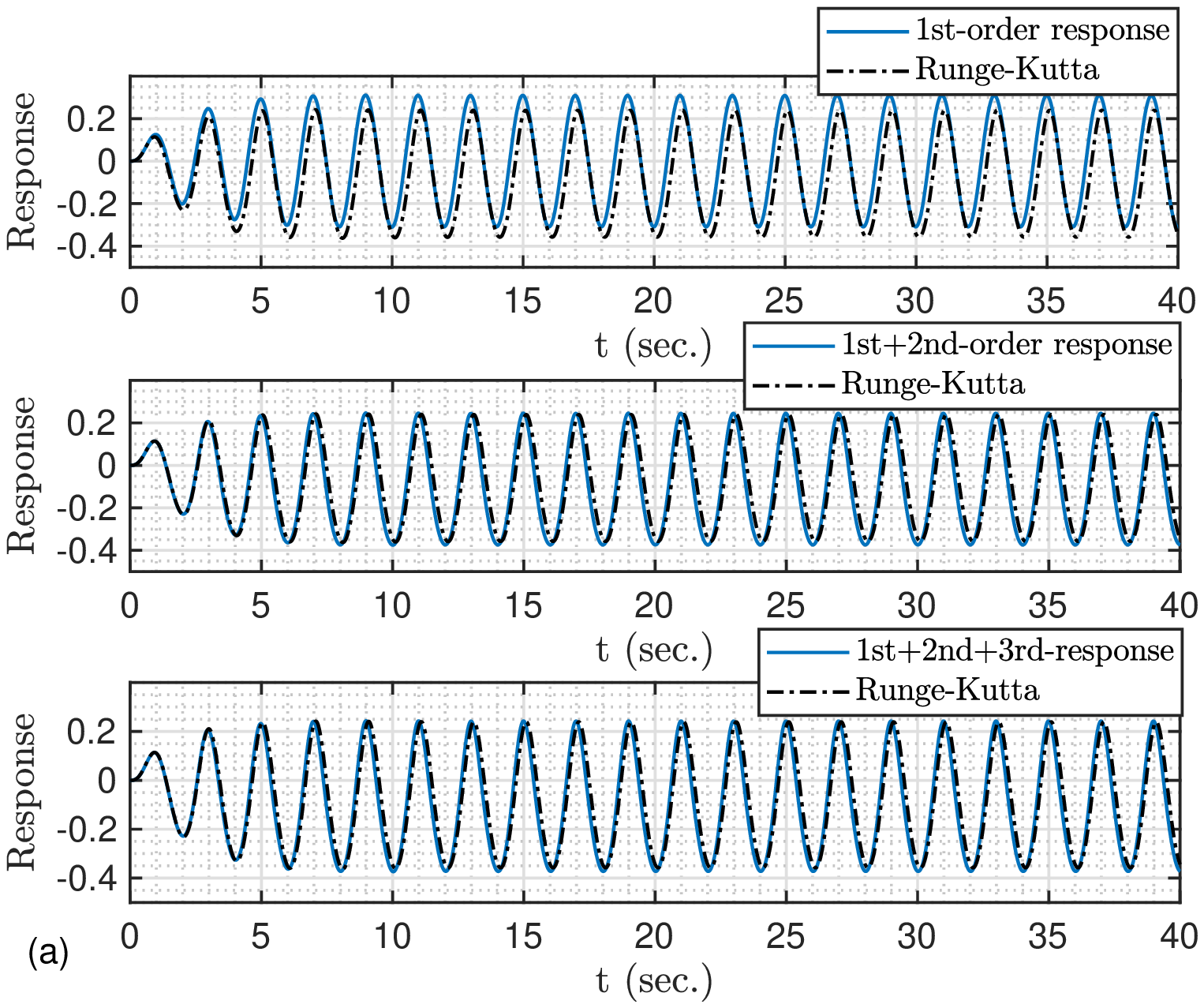}
\end{minipage}
\begin{minipage}[t]{0.5\linewidth}
\centering
\includegraphics[width=3in]{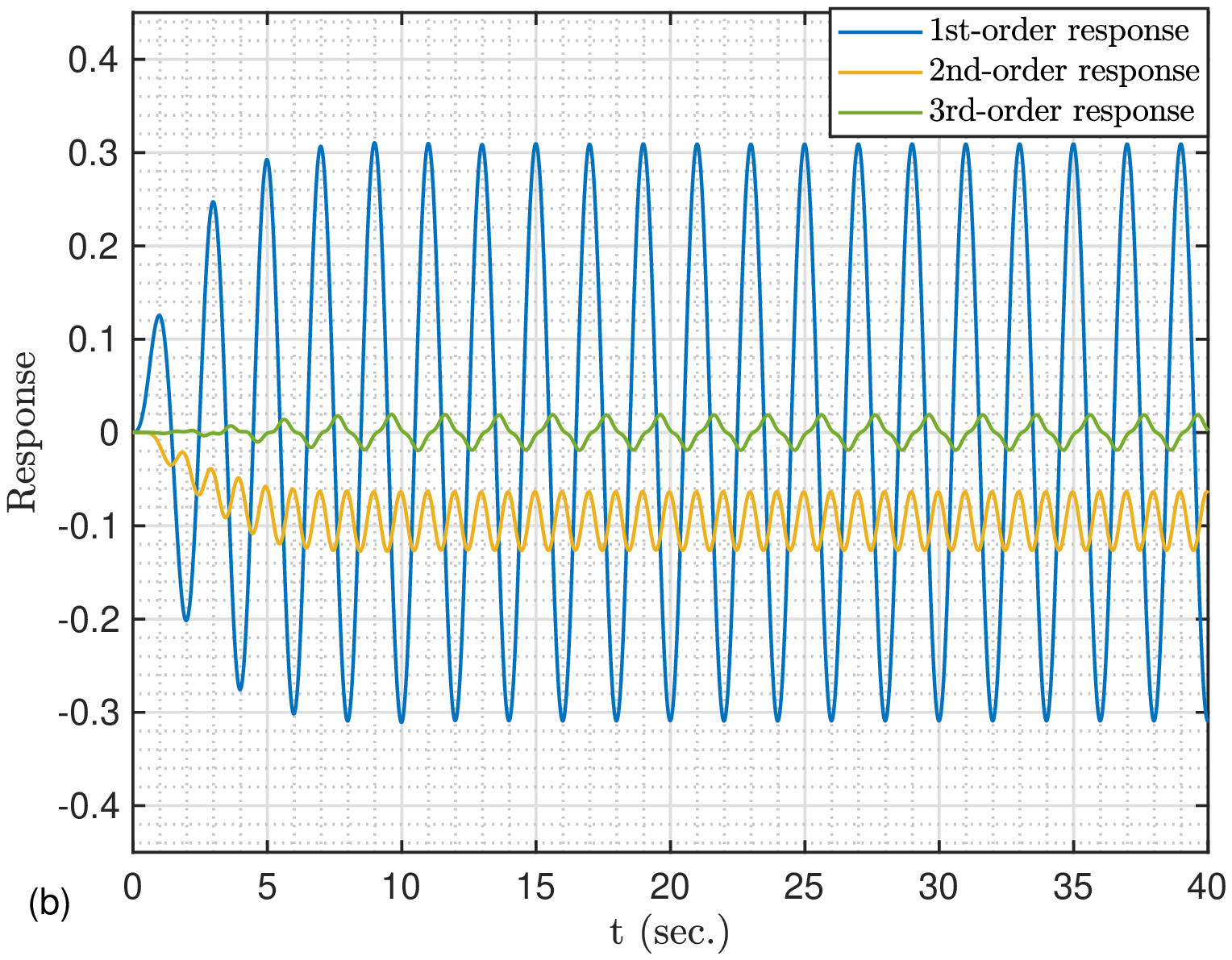}
\end{minipage}
\centering
\caption{Response for Case 3: (a) comparison between the proposed method and Runge--Kutta method, (b) contribution of the three components}\label{exp1_sinu1}
\end{figure}

\begin{figure}[H]
\centering
\begin{minipage}[t]{0.5\linewidth}
\centering
\includegraphics[width=3in]{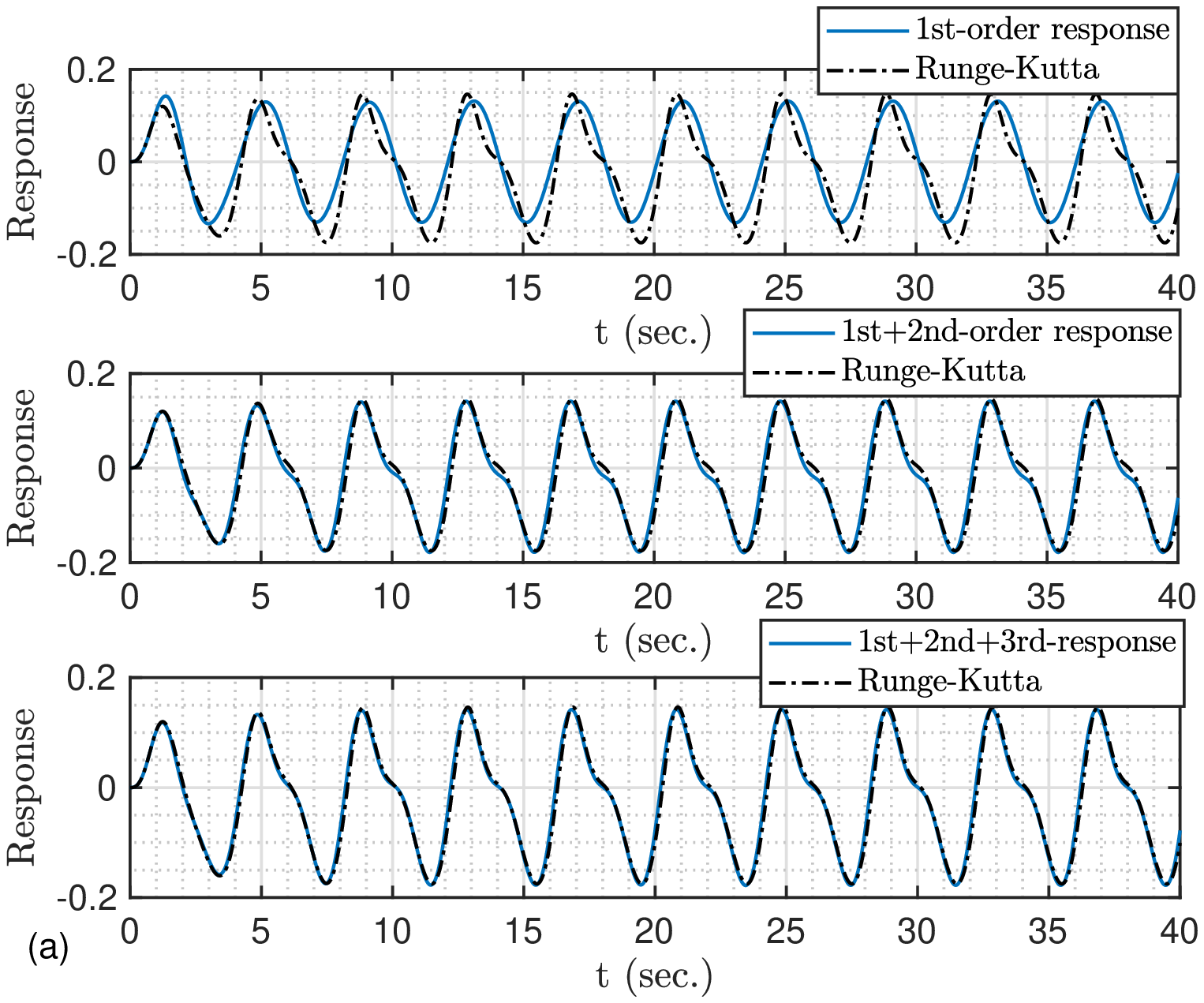}
\end{minipage}
\begin{minipage}[t]{0.5\linewidth}
\centering
\includegraphics[width=3in]{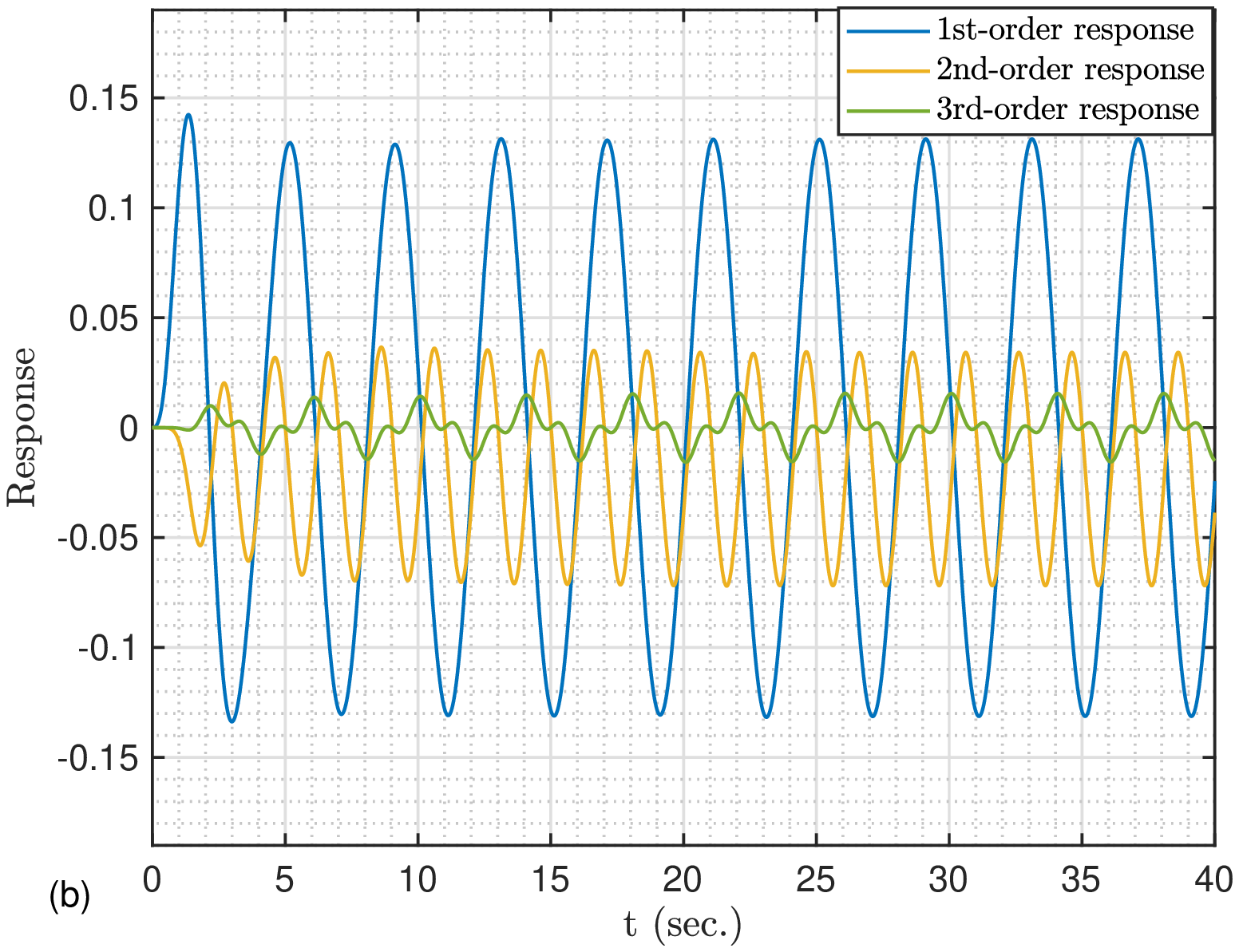}
\end{minipage}
\centering
\caption{Response for Case 4: (a) comparison between the proposed method and Runge--Kutta method, (b) contribution of the three components}\label{exp1_sinu05}
\end{figure}

\begin{figure}[H]
\centering
\begin{minipage}[t]{0.5\linewidth}
\centering
\includegraphics[width=3in]{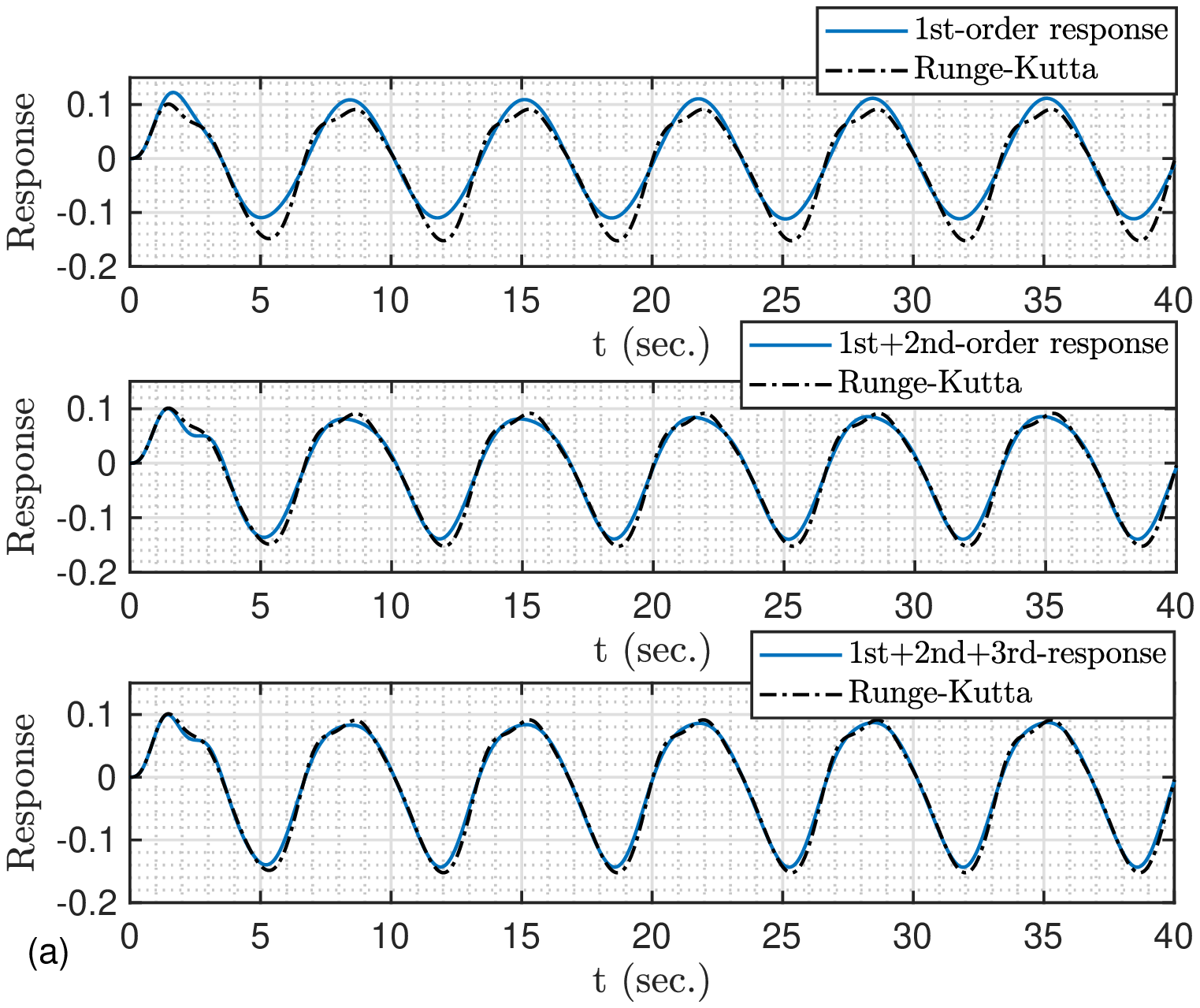}
\end{minipage}
\begin{minipage}[t]{0.5\linewidth}
\centering
\includegraphics[width=3in]{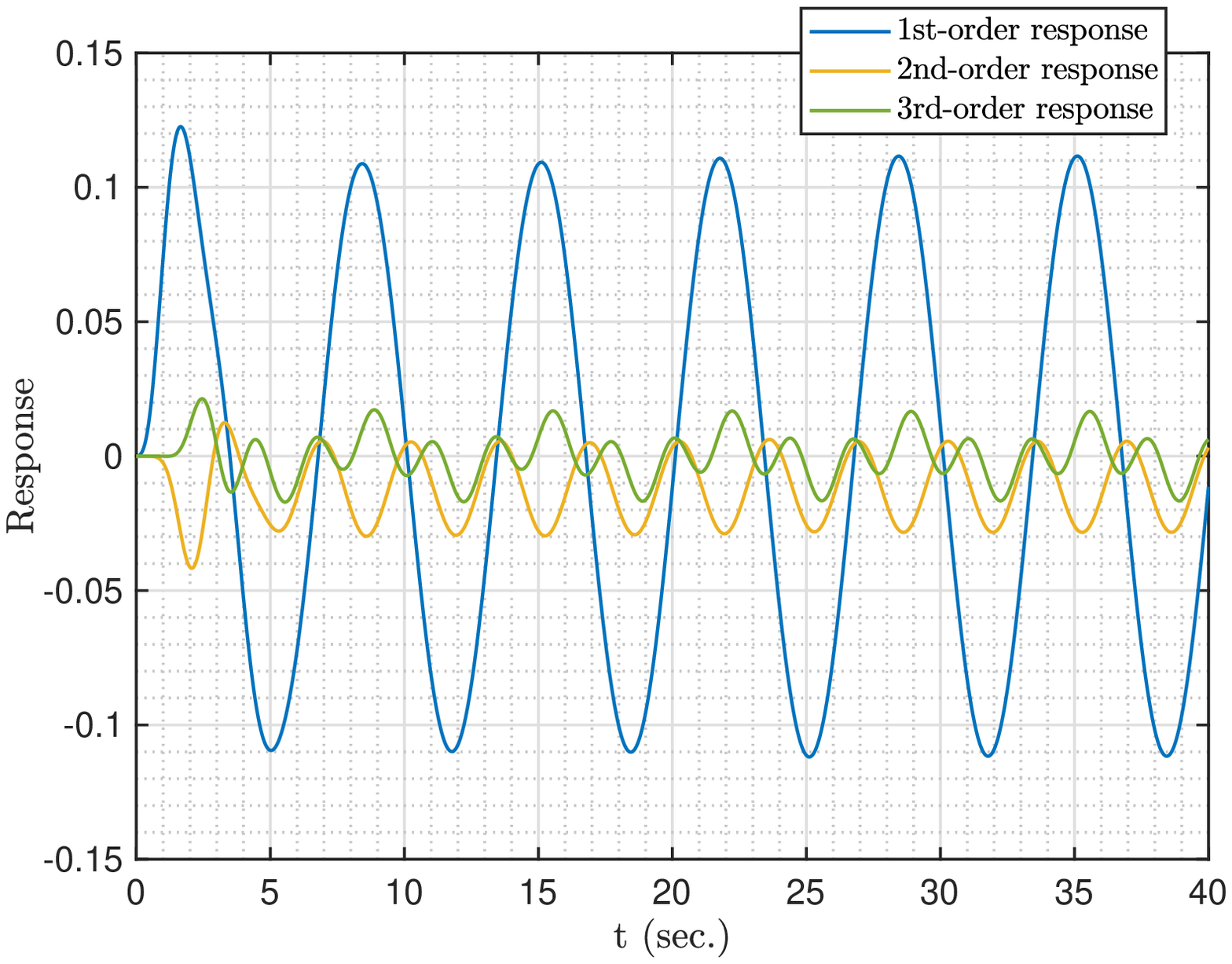}
\end{minipage}
\centering
\caption{Response for Case 5: (a) comparison between the proposed method and Runge--Kutta method, (b) contribution of the three components}\label{exp1_sinu03}
\end{figure}

Substituting poles and residues of the excitation, as well as those of the system into Eqs.~\ref{eq:residue_first} and \ref{eq:residue_last}, response coefficients $\beta_{p_i,k}$ corresponding to system poles $-a_i$ and response coefficients $\gamma_{p_i,\ell}$ corresponding to excitation poles $\lambda_\ell$ are calculated, respectively. According to Eq.~\ref{y4_final}, the first three orders of responses for each case in Table \ref{input parameters} are calculated. Figures~\ref{exp1_sinu5}(a)-\ref{exp1_sinu03}(a) show the comparison of responses obtained by the proposed method and the fourth-order Runge--Kutta method with $\Delta t=10^{-4}$. For Cases 1 and 2, the first-order responses agree well with the total responses obtained by the Runge--Kutta method, and the higher-order responses only slightly improve the transient parts. For Cases 3--5, the sum of the first three orders of responses is in good agreement with the Runge--Kutta solution. When the response nonlinearity increases, higher-order responses need to be considered. In other words, the proposed method can accurately compute the nonlinear responses by choosing a small number $N$ of Volterra series terms.

\begin{figure}[H]
\centering
\begin{minipage}[t]{0.3\linewidth}
\centering
\includegraphics[width=2in]{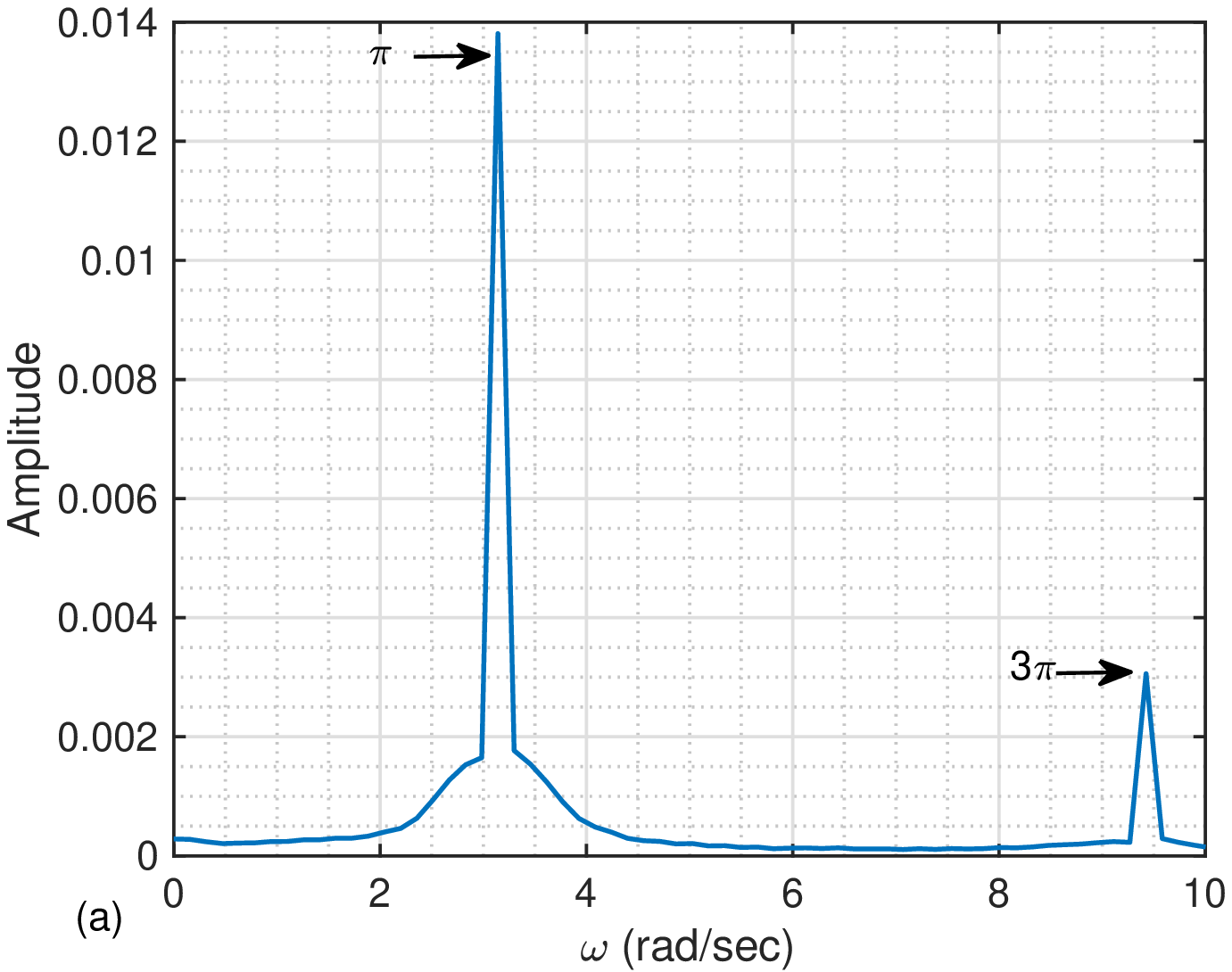}
\end{minipage}
\begin{minipage}[t]{0.3\linewidth}
\centering
\includegraphics[width=2in]{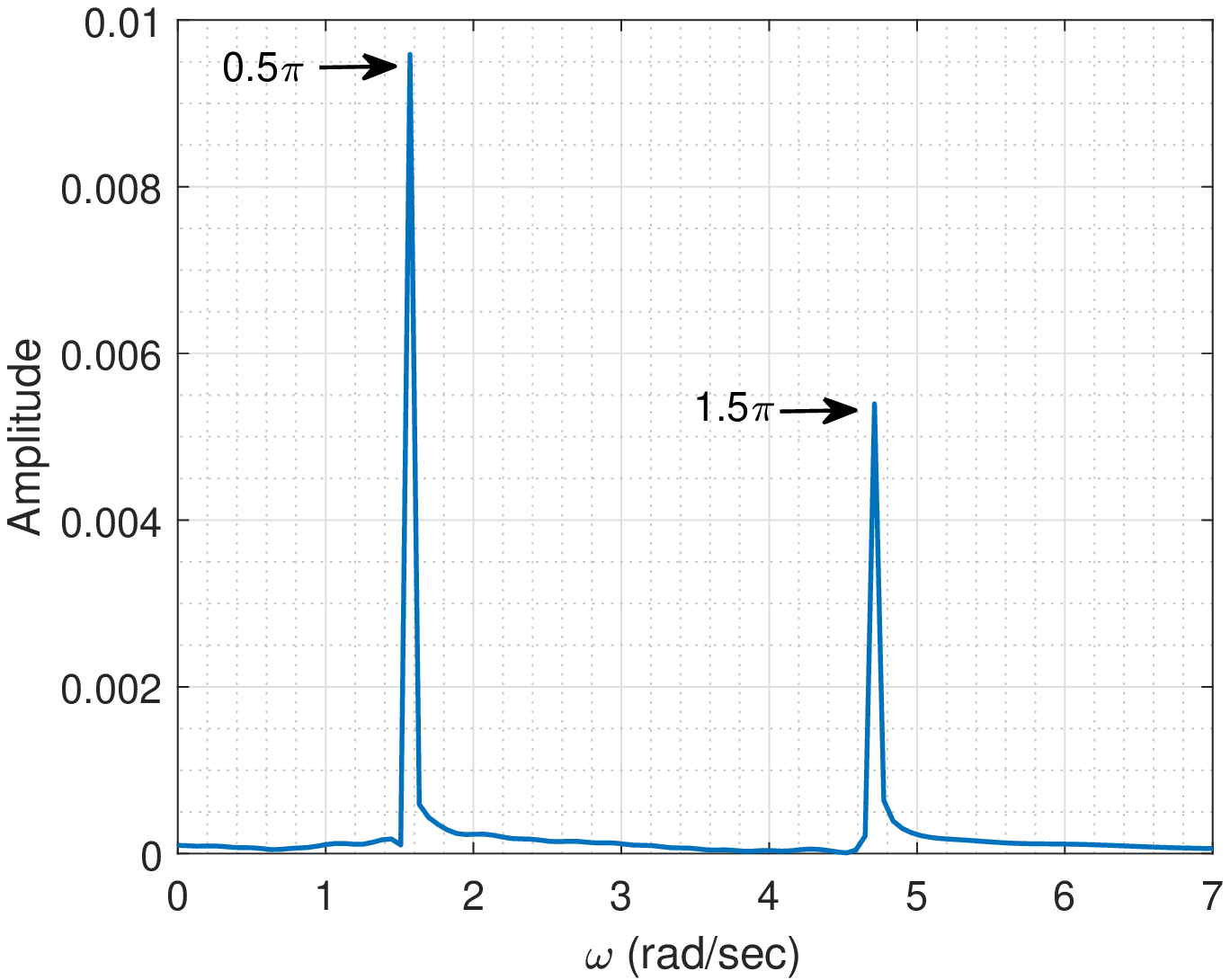}
\end{minipage}
\centering
\begin{minipage}[t]{0.3\linewidth}
\centering
\includegraphics[width=2in]{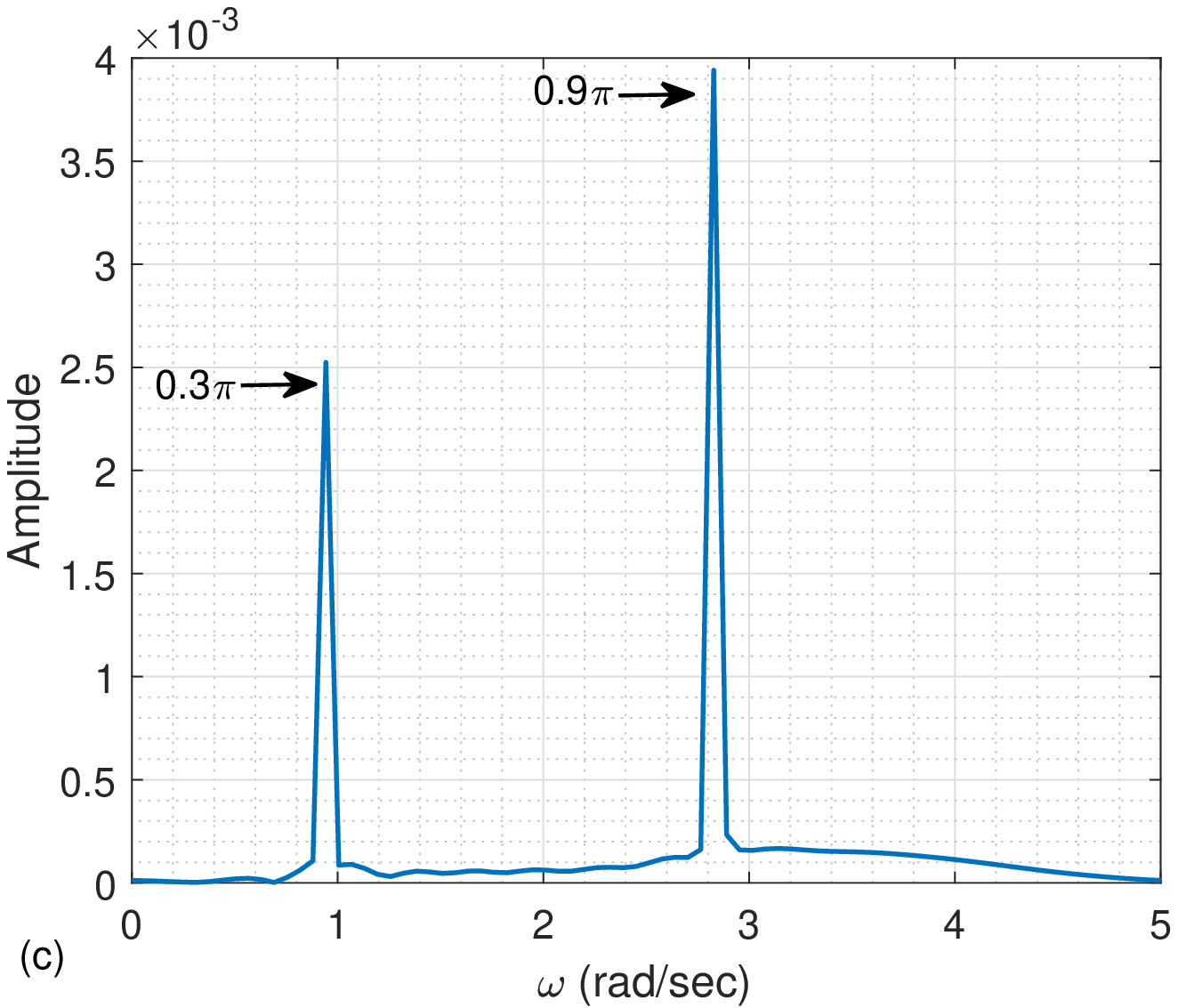}
\end{minipage}
\centering
\caption{Frequency spectra of third-order responses: (a) Case 3, (b) Case 4 and (c) Case 5}\label{Freqency_spectrum}
\end{figure}

Figures~\ref{exp1_sinu5}(b)-\ref{exp1_sinu03}(b) show the contributions of the three response components for the five cases. In each case, the first-order response is the most dominant component, and the contributions of second- and third-order responses are much less than those of the first-order response. Especially for Cases 1 and 2, whose excitation frequencies are far from the linear natural frequency, second- and third-order responses are close to zero. This may be because the QFRF and CFRF approach zero when the frequency is larger than 4 rad/s (see Figs.~\ref{QFRF}--\ref{CFRF_diff}). Furthermore, the mean values of the first-order responses are approximately zero, and those of the second-order responses are always smaller than zero, which are the difference frequency components in Eq.~\ref{y51_forced3}. Moreover, it is clearly observed that second-order responses for Cases 3--5 exhibit a periodic oscillation with a period near half of that for the first-order response, which is excited by the sum frequency component of the excitation (see second part of Eq.~\ref{y51_forced3}).
Compared with steady-state solutions of first- and second-order responses, those of third-order responses in Cases 3--5 are no longer single regular motions. By performing the FFT, frequency spectra of these three third-order responses are shown in Fig.~\ref{Freqency_spectrum}. We find that these three third-order responses are all dominated by their own fundamental harmonic component and the third harmonic (triple frequency) component.

\begin{figure}[H]
\centering
\includegraphics[width=4 in]{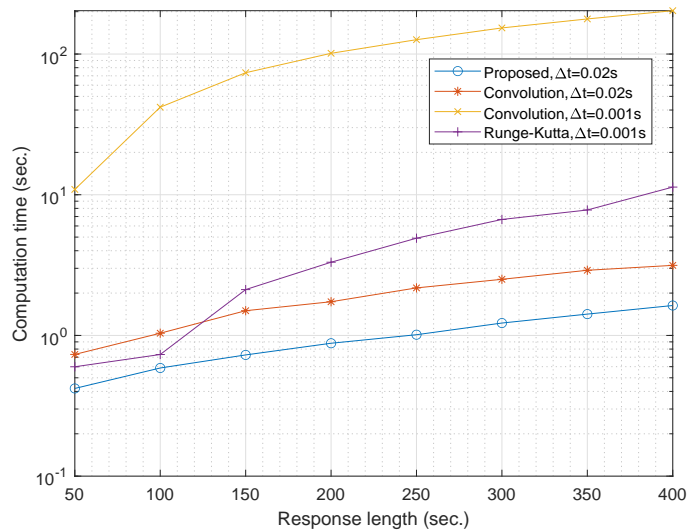}\\
\caption{Comparison of computation efficiency of the proposed method, the fourth--fifth order Runge--Kutta method and the convolution method for regular loading in Case 1}\label{Computation_time_sinu}
\end{figure}

Figure~\ref{Computation_time_sinu} shows the computational time to calculate the response of the oscillator for Case 1 by the proposed method, the fourth-order Runge--Kutta method and the convolution method. The proposed method, which has an explicit solution, is much more efficient in computational time than the latter two methods, which need small time steps to obtain high-precision solutions. In particular, the efficiency of the proposed method increases with the length of the response time.


\subsubsection{Irregular excitation}\label{irregular}
\begin{figure}[H]
\centering
\includegraphics[width=4 in]{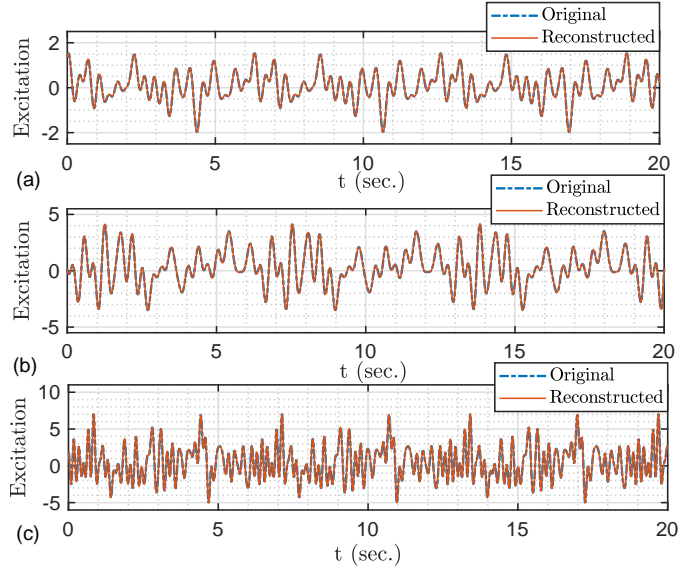}\\
\caption{Comparison of original excitations and reconstructed results: (a) Case 1, (b) Case 2 (c) Case 3}\label{exp1_whitenoise}
\end{figure}

In Eq.~\ref{eq:nonlinear_SDOF}, considering an irregular excitation consisting of several cosine functions
\begin{equation}\label{eq:irregularload}
    f(t) = \sum_{n=1}^{N_f} A_n \cos(\Omega_n t + \theta_n)
\end{equation}
where $N_f$ is the number of cosine components; $A_n$, $\Omega_n$ and $\theta_n$ are the amplitude, frequency and phase angle of the $n^{th}$ component, respectively. Table ~\ref{input parameters2} lists three cases of these parameters. In each case, the amplitudes of all components are the same, and phase angles $\theta_n$ uniformly distributed between 0 and $2 \pi$ are randomly generated. To decompose the excitation into a pole-residue form, the Prony-SS method is used, whose concept is similar to that of a principal component method. The readers are referred to Ref.~\cite{hu2013signal} for details. The chosen rank of each case is also shown in Table ~\ref{input parameters2}. Figure~\ref{exp1_whitenoise} shows the comparison of original excitations and reconstructed results of these three cases, which all have excellent agreement.
\begin{table}[H]
   \begin{center}
   \caption{Parameter values of the irregular excitation}
   \label{input parameters2}
   \begin{tabular}{|c|c|c|c|c|}
   \hline
   Case &  $A_n$ & $\Omega_n$   & $\theta_n$  & Rank \\
   \hline
   $1$      & $0.2$  &  $[0:1:20]$   & Uniform Random Number  &  $42$\\
   \hline
   $2$      & $0.5$  &  $[0:1:20]$   & Uniform Random Number  & $42$\\
   \hline
   $3$      & $0.5$  &  $[0:1:40]$   & Uniform Random Number  & $82$\\
   \hline
   \end{tabular}
   \end{center}
\end{table}

\begin{figure}[H]
\centering
\begin{minipage}[t]{0.5\linewidth}
\centering
\includegraphics[width=3in]{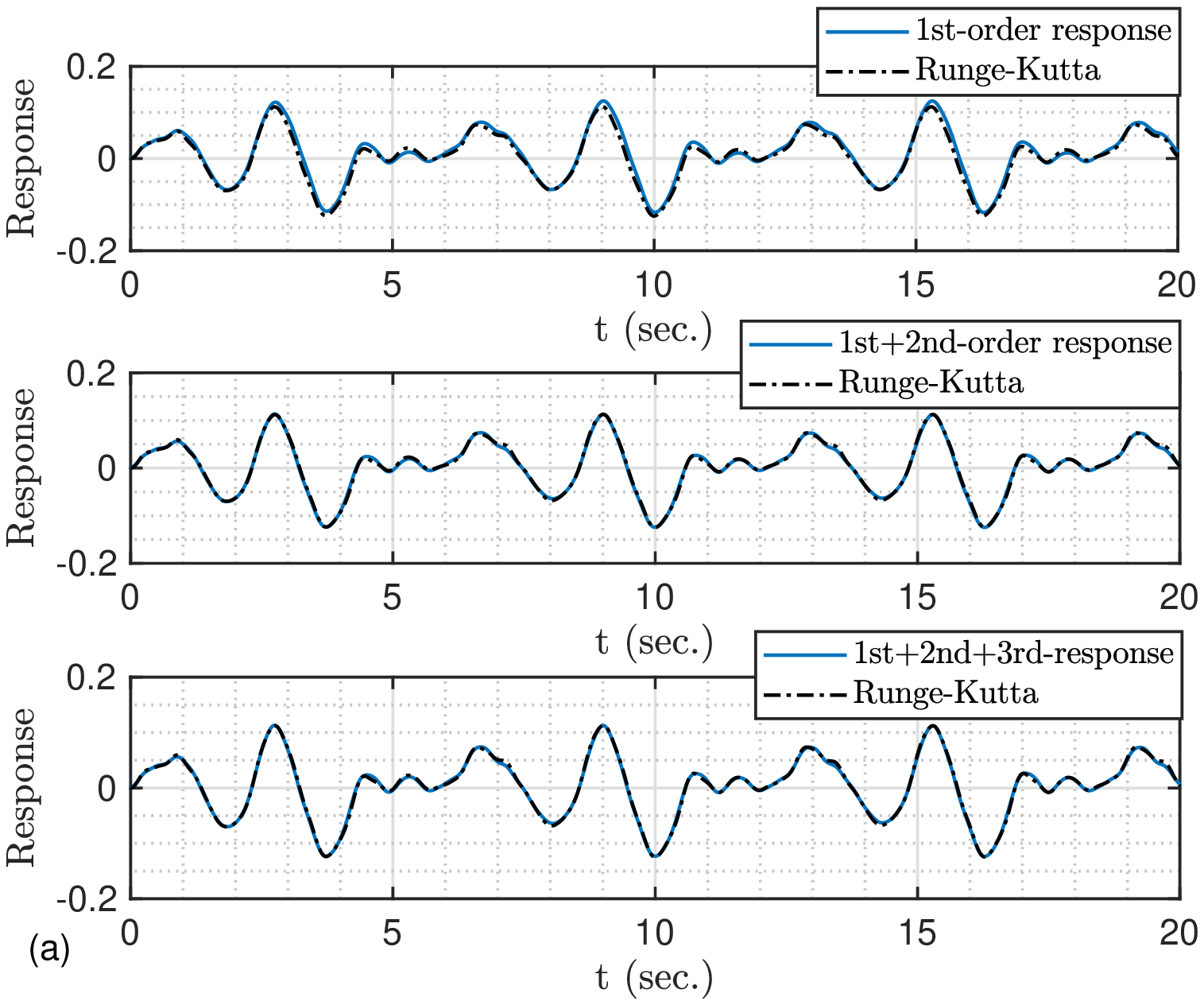}
\end{minipage}
\begin{minipage}[t]{0.5\linewidth}
\centering
\includegraphics[width=3in]{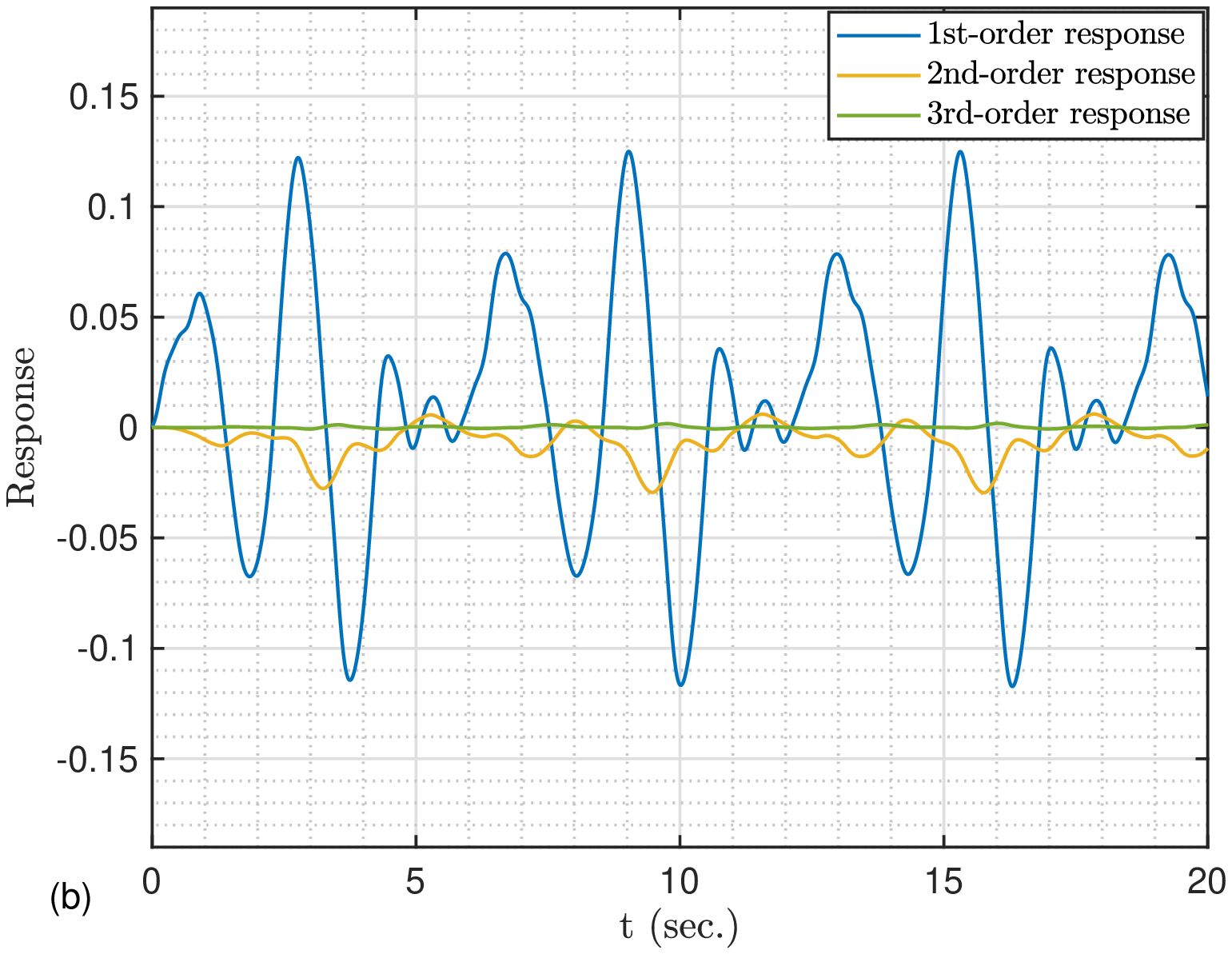}
\end{minipage}
\centering
\caption{Response to irregular excitation for Case 1: (a) comparison between the proposed method and Runge--Kutta method, (b) contribution of the three components}\label{exp1_whitenoise_A02}
\end{figure}

Referring to Eq.~\ref{y4_final} with $N=3$, nonlinear responses of the system to irregular excitations in Table ~\ref{input parameters2} are calculated, which are shown in Figs.~\ref{exp1_whitenoise_A02}(a)-\ref{exp1_whitenoise_A05}(a), respectively. Additionally, Figs.~\ref{exp1_whitenoise_A02}(a)-\ref{exp1_whitenoise_A05}(a) show the results computed by the fourth-order Runge--Kutta method. In all cases, the sums of the first three orders of responses agree well with those obtained by the Runge--Kutta method.

The contributions of the first three orders of responses for each case are plotted in Figs.~\ref{exp1_whitenoise_A02}(b)-\ref{exp1_whitenoise_A05}(b). Similarly, the system vibration is dominated by the first-order response. However, the contributions of second- and third-order responses significantly grow with increasing excitation magnitude and frequency number. Furthermore, when the magnitude of the nonlinear response becomes large, sharp troughs are present. This phenomenon may be induced by the nonlinear stiffness. While the first-order response fails to capture these troughs, the higher-order responses successfully capture these troughs.

\begin{figure}[H]
\centering
\includegraphics[width=4 in]{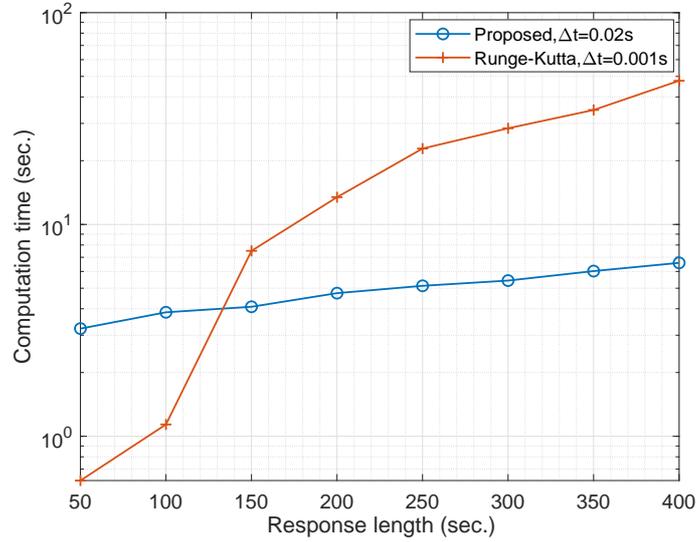}\\
\caption{Comparison of computation efficiency of the proposed method and the fourth--fifth order Runge--Kutta method for irregular loading in Case 1}\label{Computation_time_irregular}
\end{figure}

Figure~\ref{Computation_time_irregular} plots the computational time to calculate the response of the oscillator for the irregular loading in Case 1 by the proposed method and the fourth--fifth order Runge--Kutta method, respectively. While the fourth--fifth order Runge--Kutta method is more efficient under a small response length, the proposed method becomes much more efficient when the response length is larger than about 130 s. In addition, the proposed method obtains the explicit response solution, so one can directly obtain the response value at a specific time $t_p$ instead of integrating from $0$ to $t_p$ for traditional numerical methods.
\begin{figure}[H]
\centering
\begin{minipage}[t]{0.5\linewidth}
\centering
\includegraphics[width=3in]{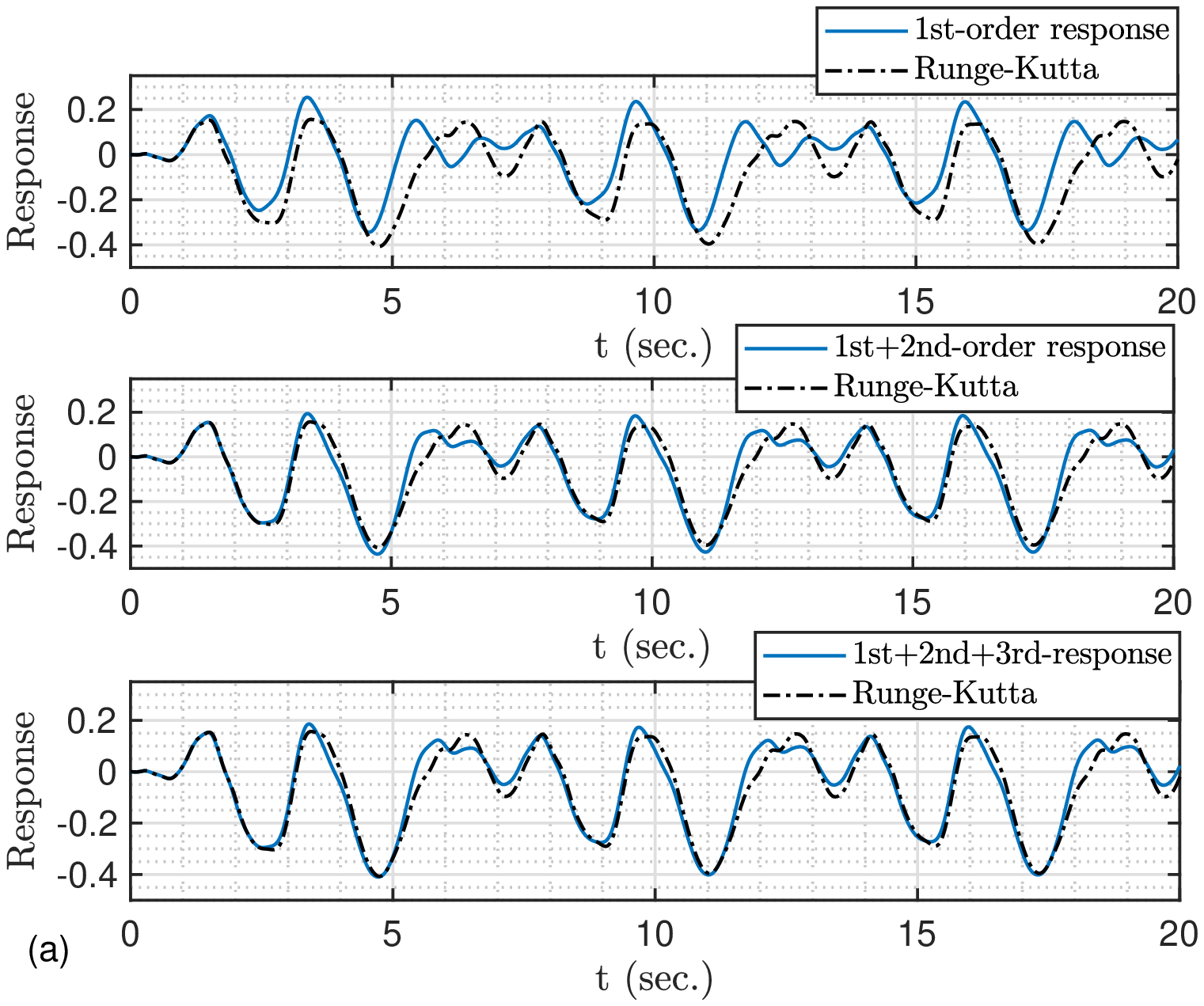}
\end{minipage}
\begin{minipage}[t]{0.5\linewidth}
\centering
\includegraphics[width=3in]{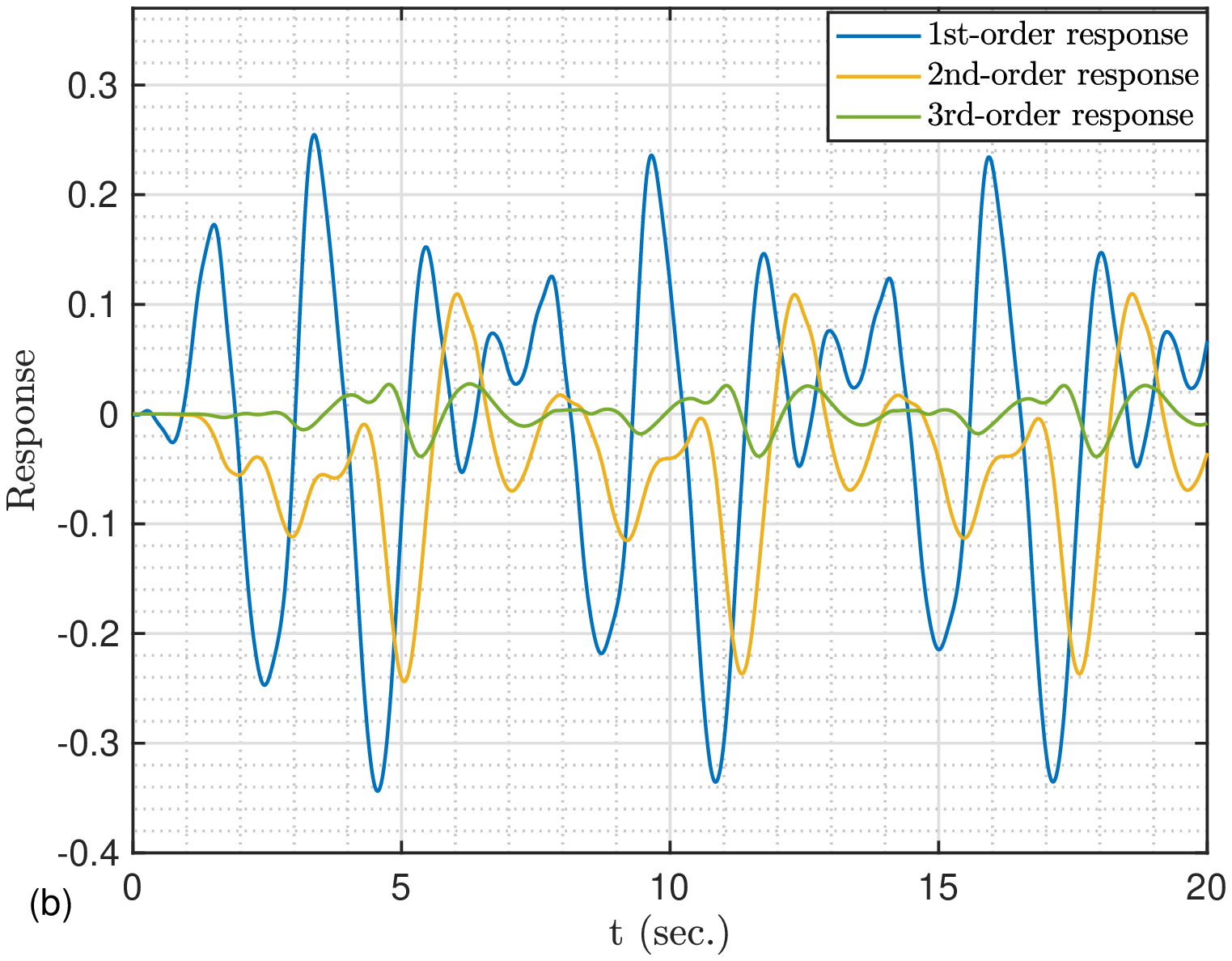}
\end{minipage}
\centering
\caption{Response to irregular excitation for Case 2: (a) comparison between the proposed method and Runge--Kutta method, (b) contribution of the three components}\label{exp1_whitenoise_A05}
\end{figure}
\begin{figure}[H]
\centering
\begin{minipage}[t]{0.5\linewidth}
\centering
\includegraphics[width=3in]{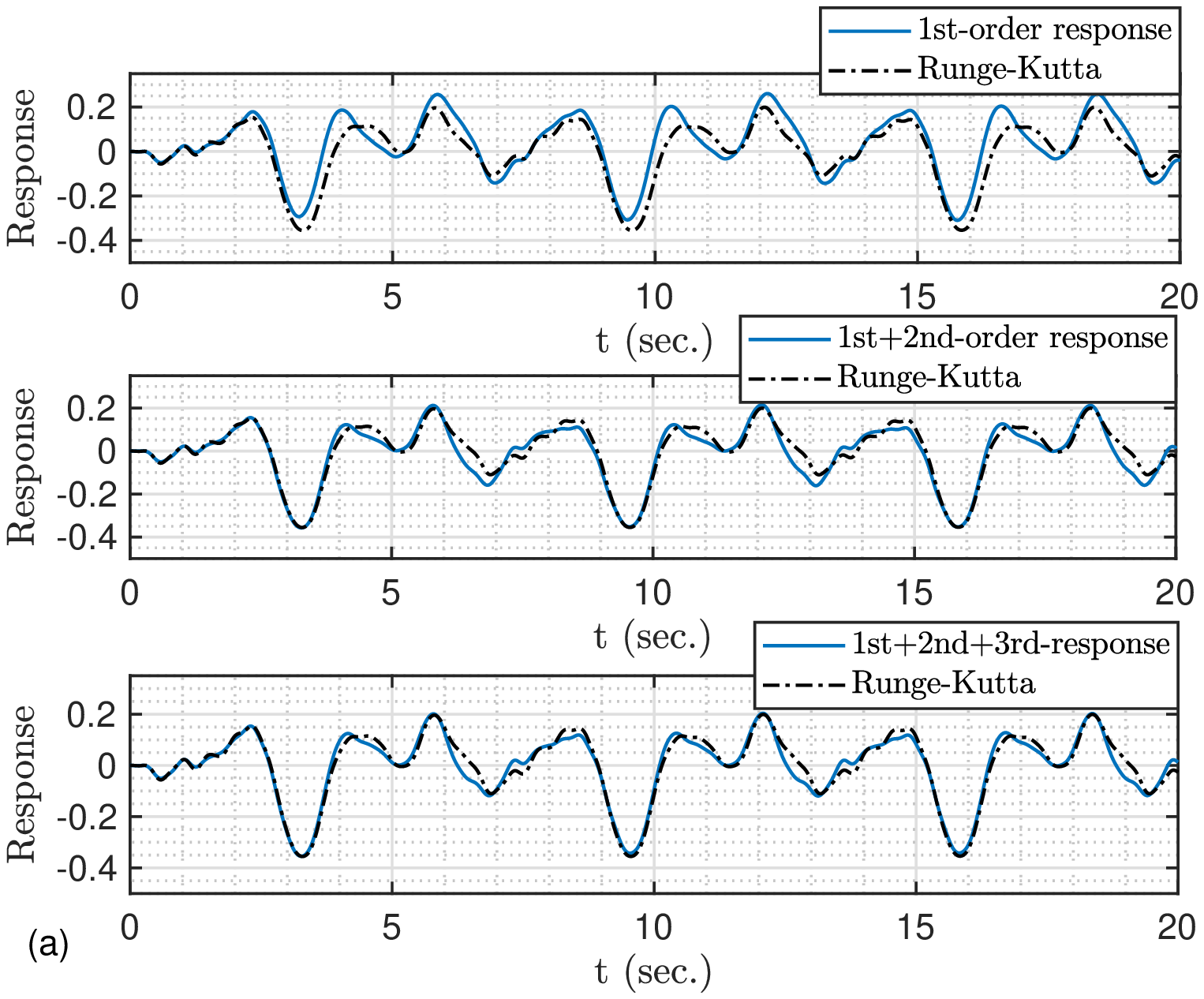}
\end{minipage}
\begin{minipage}[t]{0.5\linewidth}
\centering
\includegraphics[width=3in]{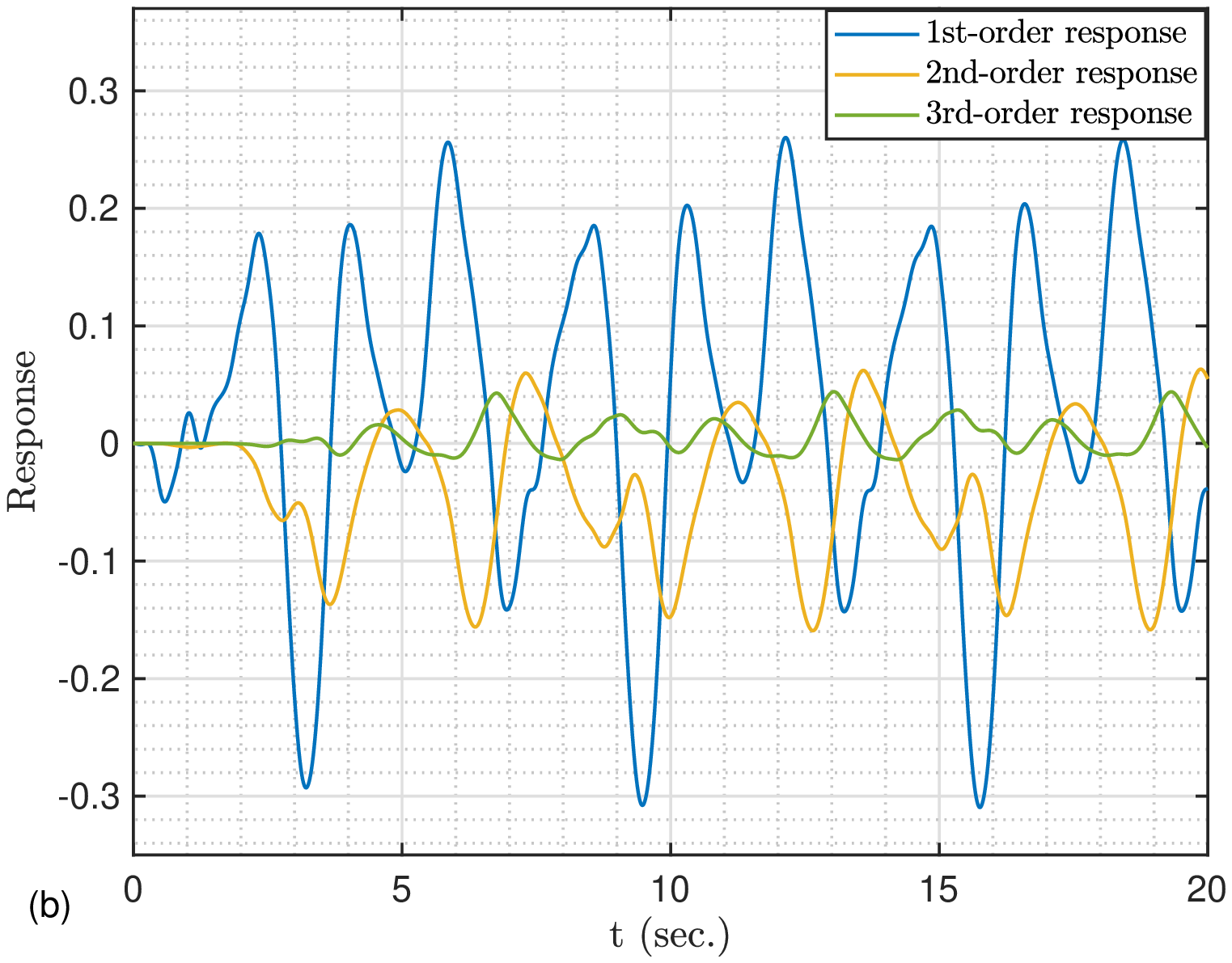}
\end{minipage}
\centering
\caption{Response to irregular excitation for Case 3: (a) comparison between the proposed method and Runge--Kutta method, (b) contribution of the three components}\label{exp1_whitenoise_A05}
\end{figure}

\subsection{An unknown nonlinear system}
\begin{figure}[H]
\centering
\includegraphics[width=4in]{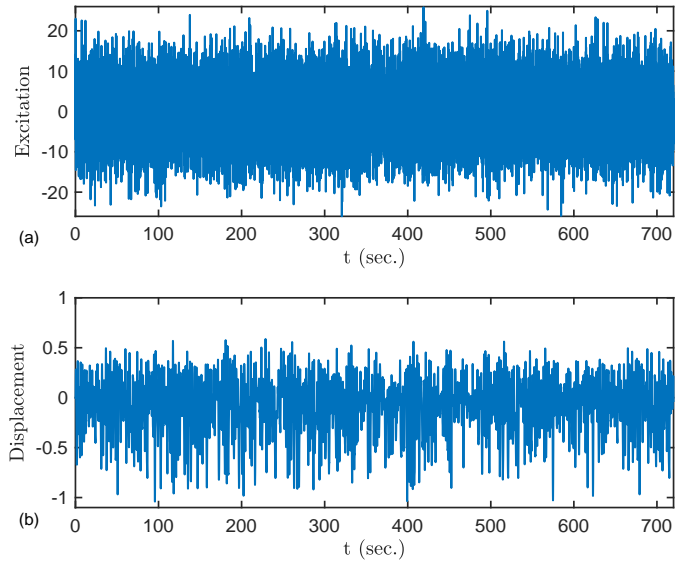}
\caption{Input--output dataset used to identify Volterra series: (a) input excitation, (b) output response}\label{Exam2_input1}
\end{figure}

\begin{figure}[H]
\centering
\begin{minipage}[t]{0.5\linewidth}
\centering
\includegraphics[width=3in]{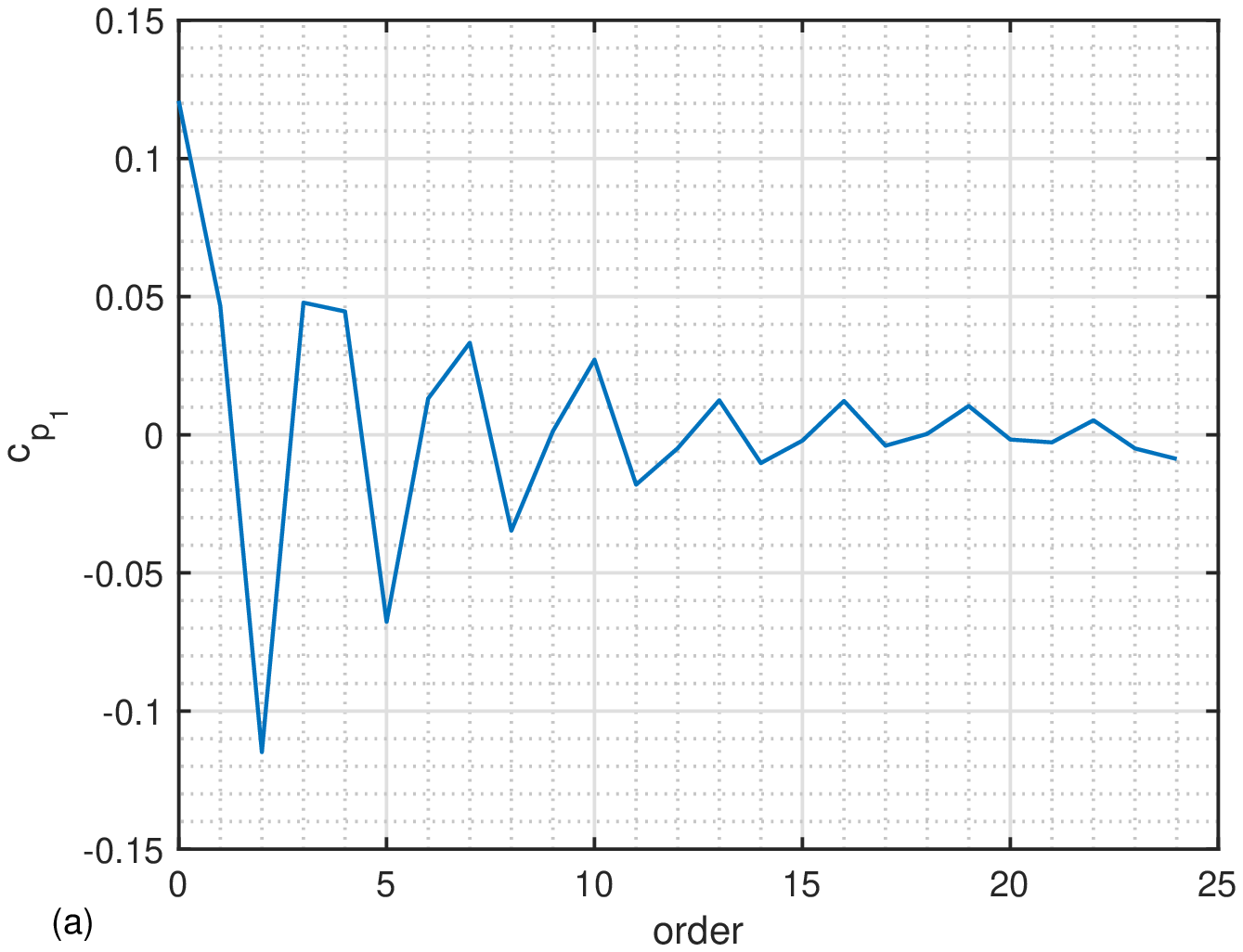}
\end{minipage}
\begin{minipage}[t]{0.5\linewidth}
\centering
\includegraphics[width=3in]{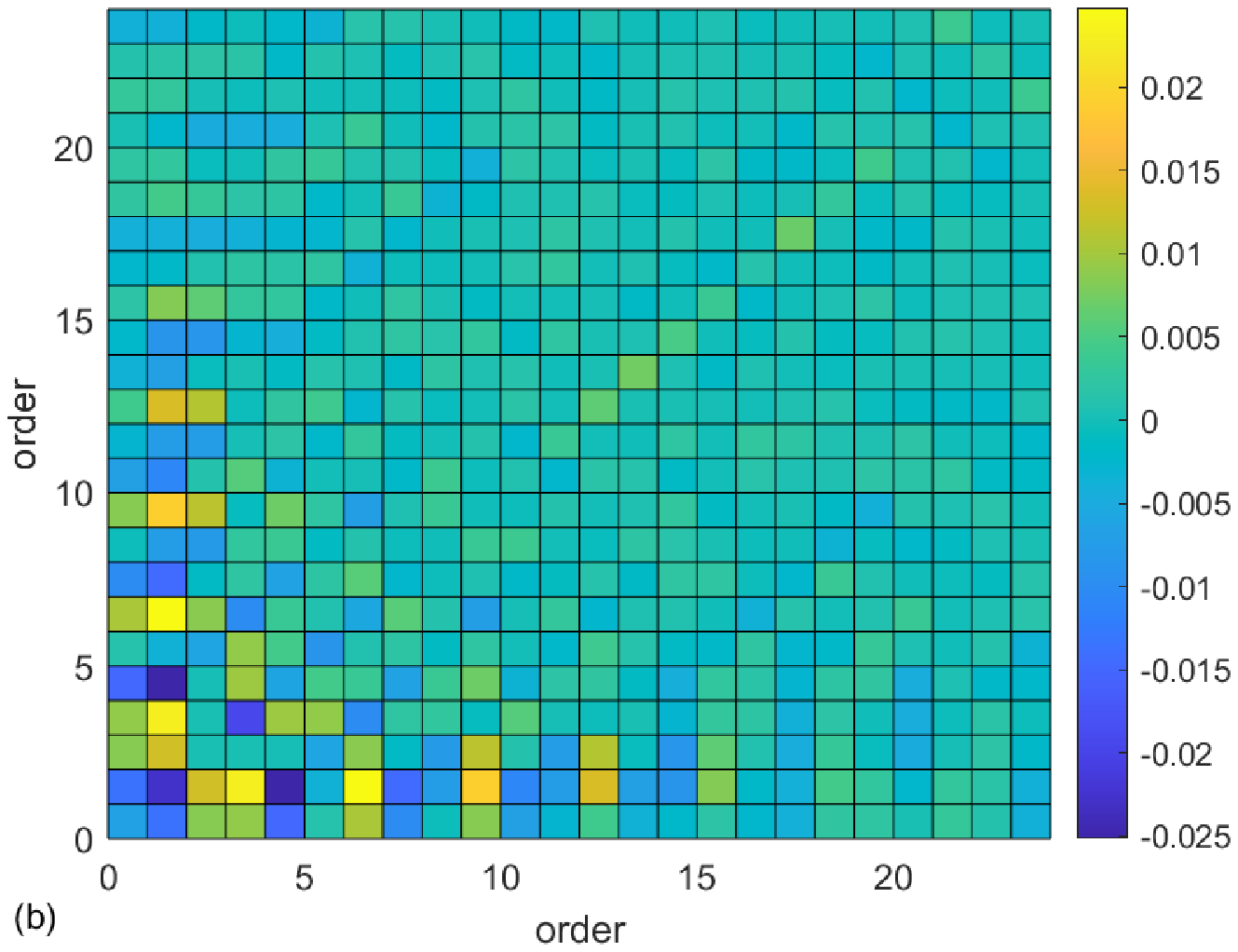}
\end{minipage}
\centering
\caption{Coefficients of the Volterra series: (a) $c_{p_1}$, (b) $c_{p_1p_2}$}\label{Exam2_coef}
\end{figure}

\begin{figure}[H]
\centering
\begin{minipage}[t]{0.5\linewidth}
\centering
\includegraphics[width=3in]{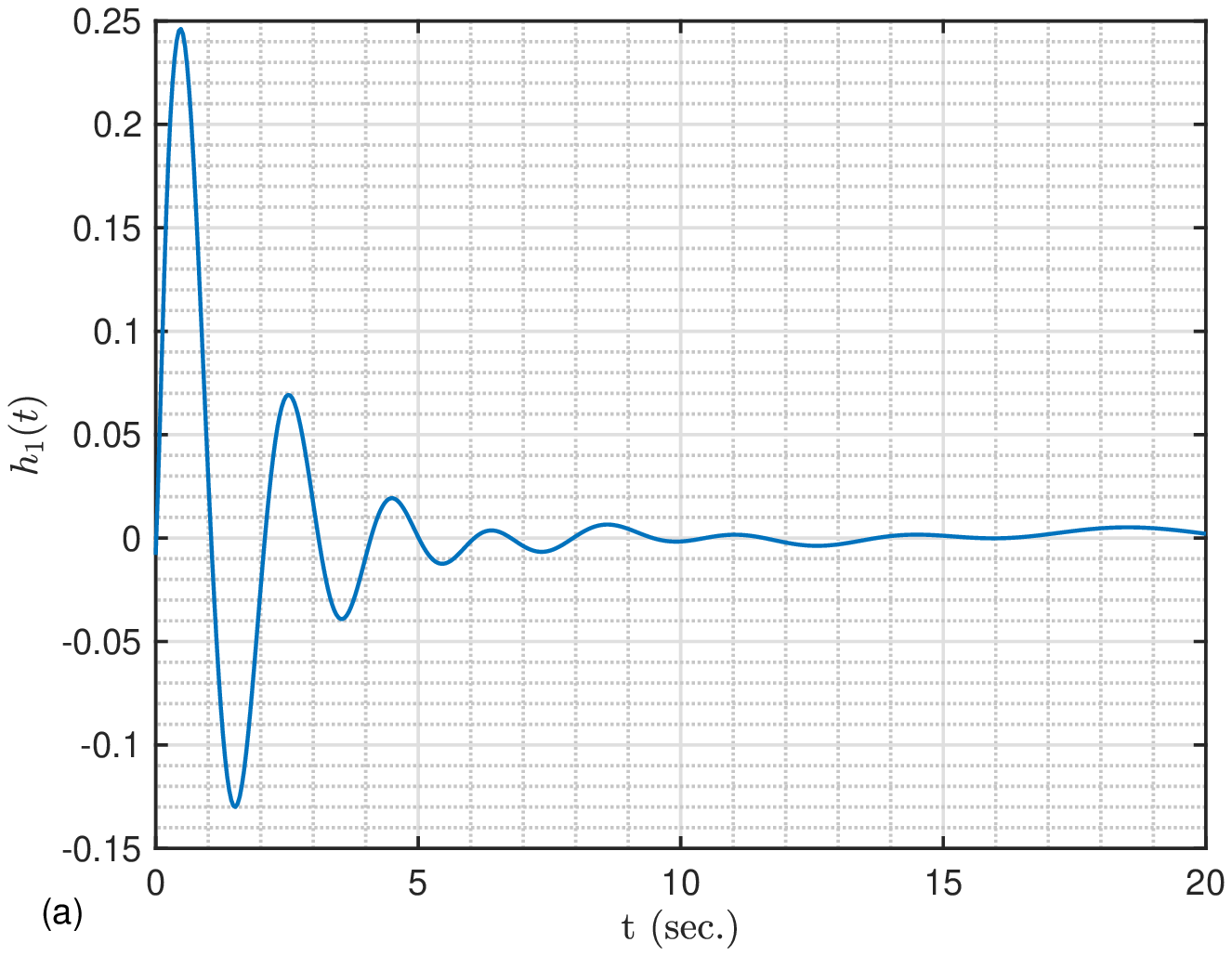}
\end{minipage}
\begin{minipage}[t]{0.5\linewidth}
\centering
\includegraphics[width=3in]{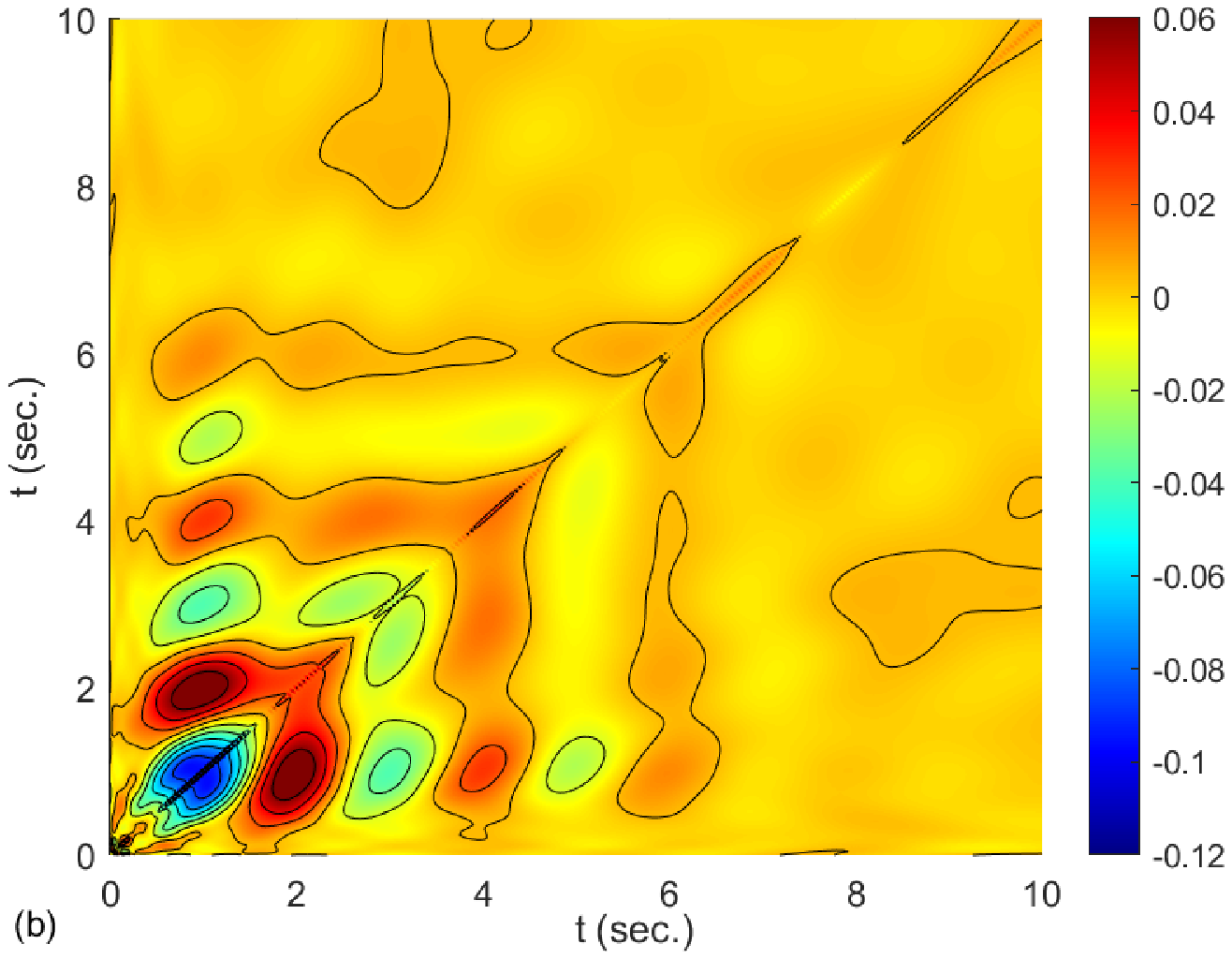}
\end{minipage}
\centering
\caption{Identified Volterra series: (a) $h_1(t)$, (b) $h_2(t_1, t_2)$}\label{Exam2_IRF}
\end{figure}

To check the applicability of the proposed method to an unknown nonlinear system, a known input excitation and its corresponding response are used to identify its Volterra kernel functions. When the Volterra kernel functions are known, we can follow the procedure in Section ~\ref{example1} to predict system responses. In this study, the input excitation is white noise with a constant power spectrum $S_0=0.001$, and the corresponding response is obtained by solving Eq.~\ref{eq:nonlinear_SDOF} by the fourth-order Runge--Kutta method, which is shown in Fig.~\ref{Exam2_input1}. From Section ~\ref{example1}, we determine that the sum of the first two orders of responses agrees well with the total response. In this study, the order of Volterra series $N$ is chosen to be 2, damping rates of Laguerre polynomials are $a_1=a_2=2$, and numbers of Laguerre polynomials are $R_1=R_2=24$. To estimate the first two orders of Volterra kernel functions, a matrix equation is constructed using excitation data and response data. By using the least square method~\cite{son2020parametric} to solve this matrix equation, coefficients $c_{p_1}$ and $c_{p_1p_2}$ in Eq.~\ref{y2} are identified. Figure~\ref{Exam2_coef} plots $c_{p_1}$ and $c_{p_1p_2}$, respectively, which have good agreement with the exact results in Fig.~\ref{coef_exp1}. Then, the first two order Volterra kernel functions are constructed by Eq.~\ref{ht_exp}. Compared with the exact results in Figs.~\ref{LIRF} and \ref{QIRF}, the identified Volterra kernel functions in Fig.~\ref{Exam2_IRF} completely agree well with the exact solutions. Note that the white noise excitation, which can excite more frequency components of the response, is chosen to obtain good Volterra kernel functions.

\begin{figure}[H]
\centering
\begin{minipage}[t]{0.5\linewidth}
\centering
\includegraphics[width=3in]{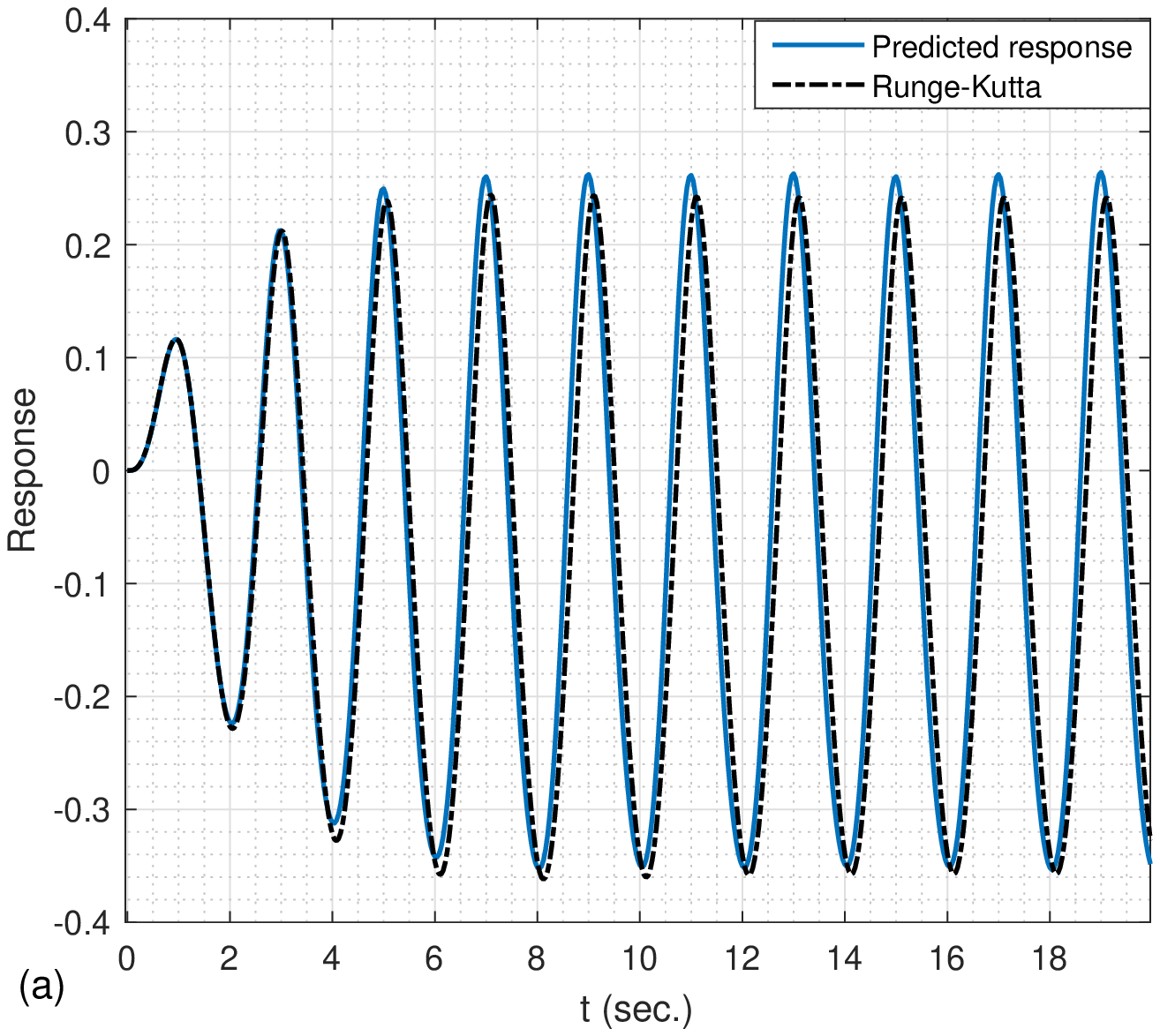}
\end{minipage}
\begin{minipage}[t]{0.5\linewidth}
\centering
\includegraphics[width=3in]{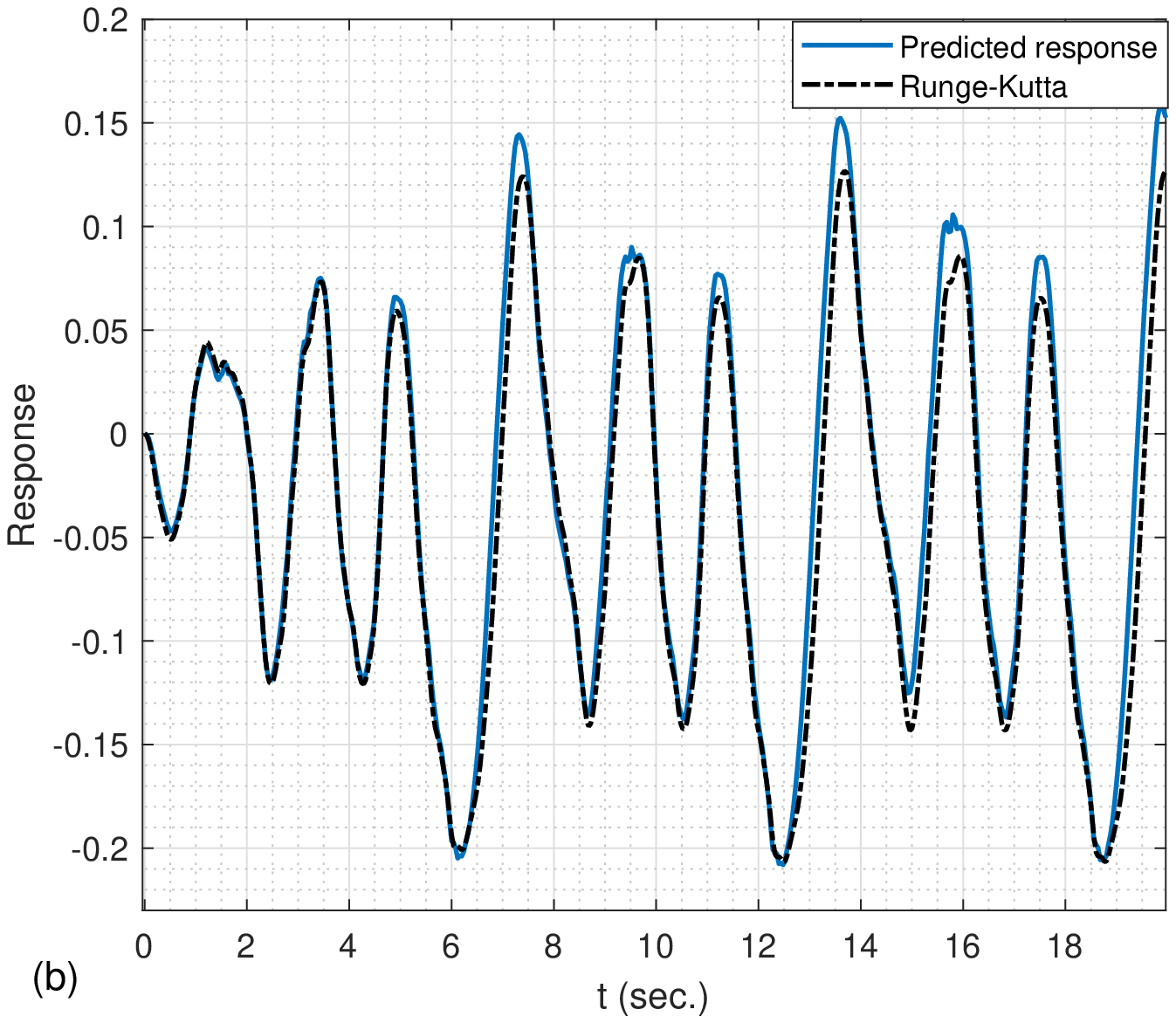}
\end{minipage}
\centering
\caption{Comparison of responses between the predicted and numerical results: (a) response to regular excitation, (b) response to irregular excitation}\label{Exam2_output23}
\end{figure}

\begin{figure}[H]
\centering
\begin{minipage}[t]{0.5\linewidth}
\centering
\includegraphics[width=3in]{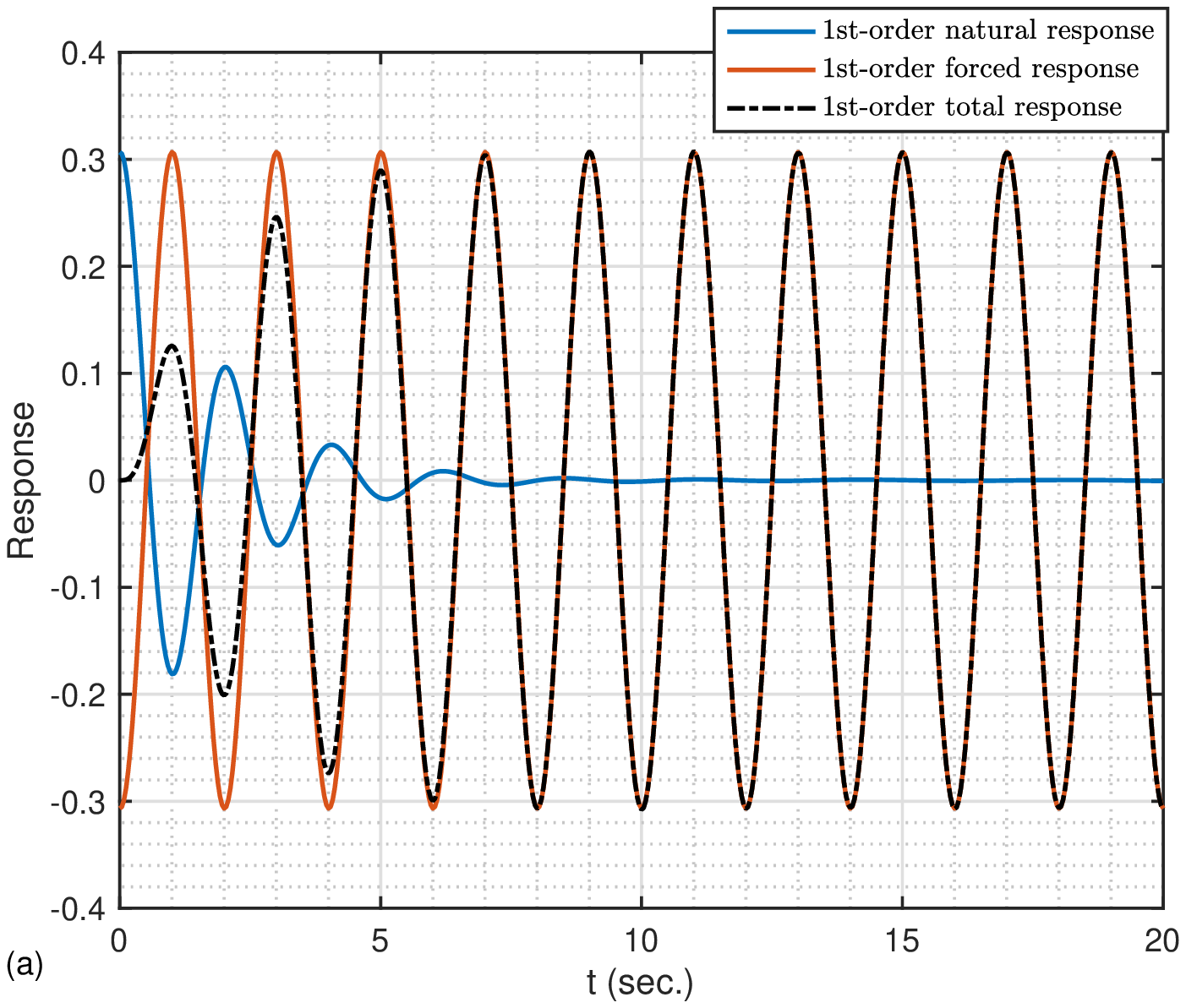}
\end{minipage}
\begin{minipage}[t]{0.5\linewidth}
\centering
\includegraphics[width=3in]{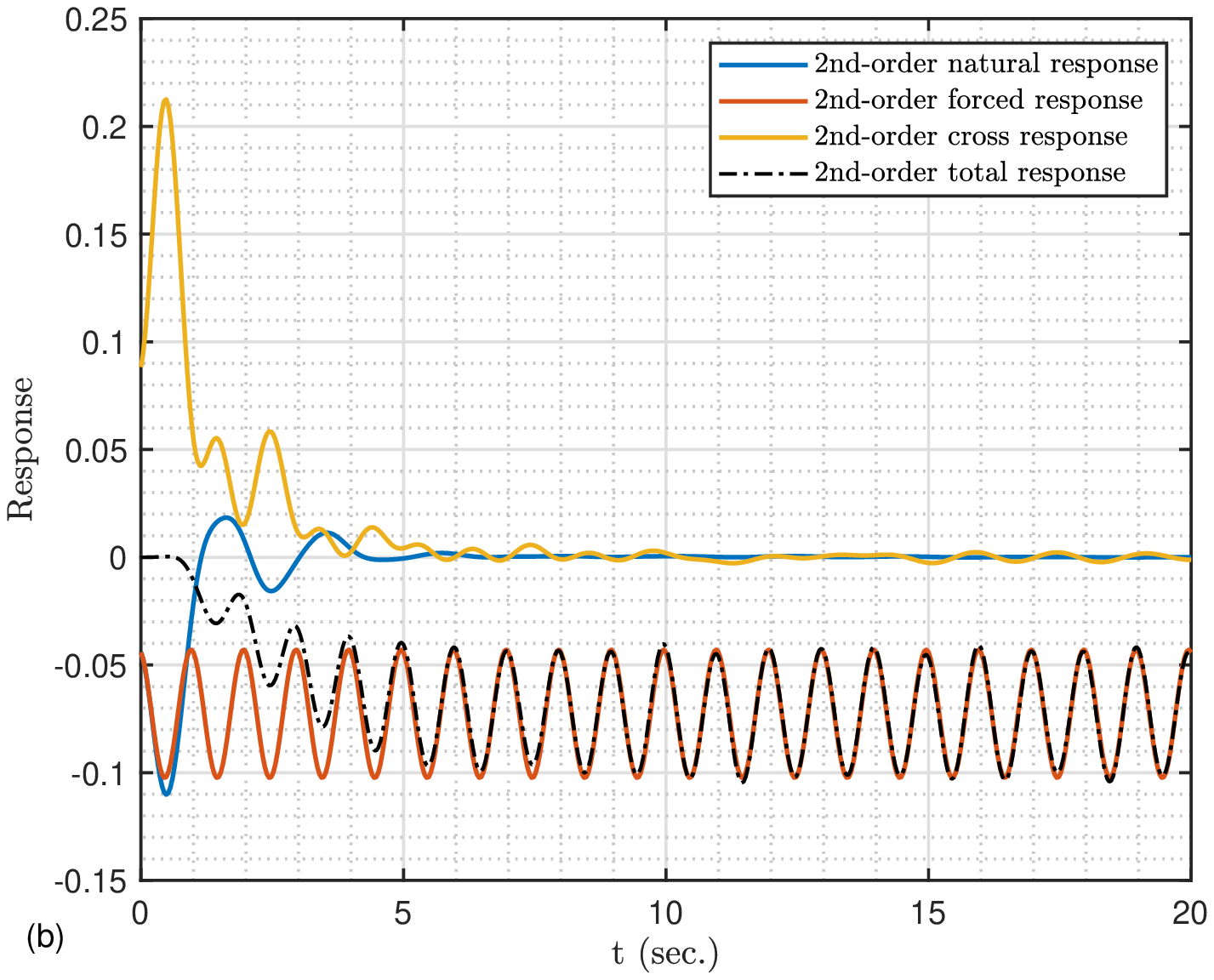}
\end{minipage}
\centering
\caption{Three components and the total responses: (a) first order and (b) second order}\label{Exam2_output2_component}
\end{figure}

\begin{figure}[H]
\centering
\begin{minipage}[t]{0.5\linewidth}
\centering
\includegraphics[width=3in]{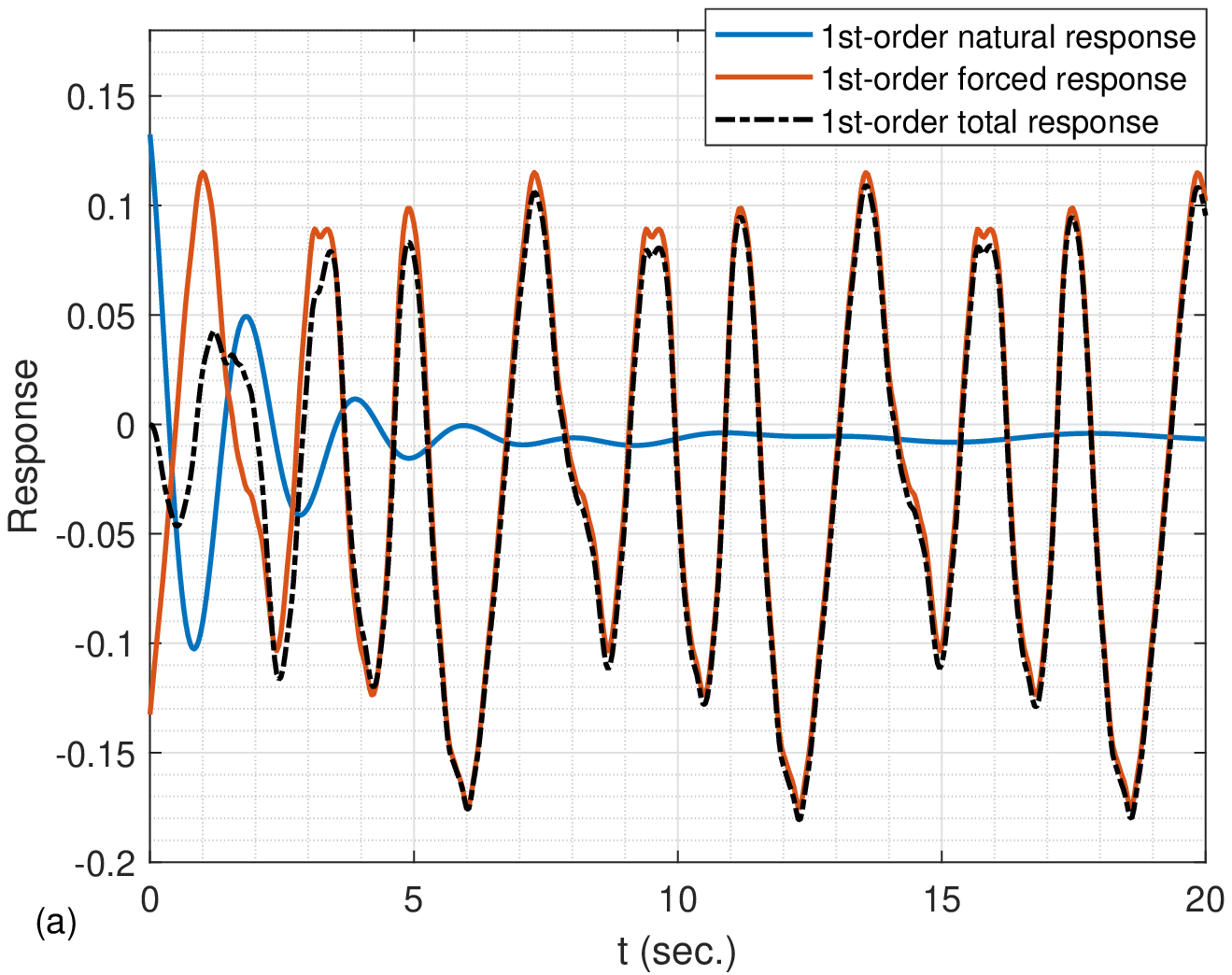}
\end{minipage}
\begin{minipage}[t]{0.5\linewidth}
\centering
\includegraphics[width=3in]{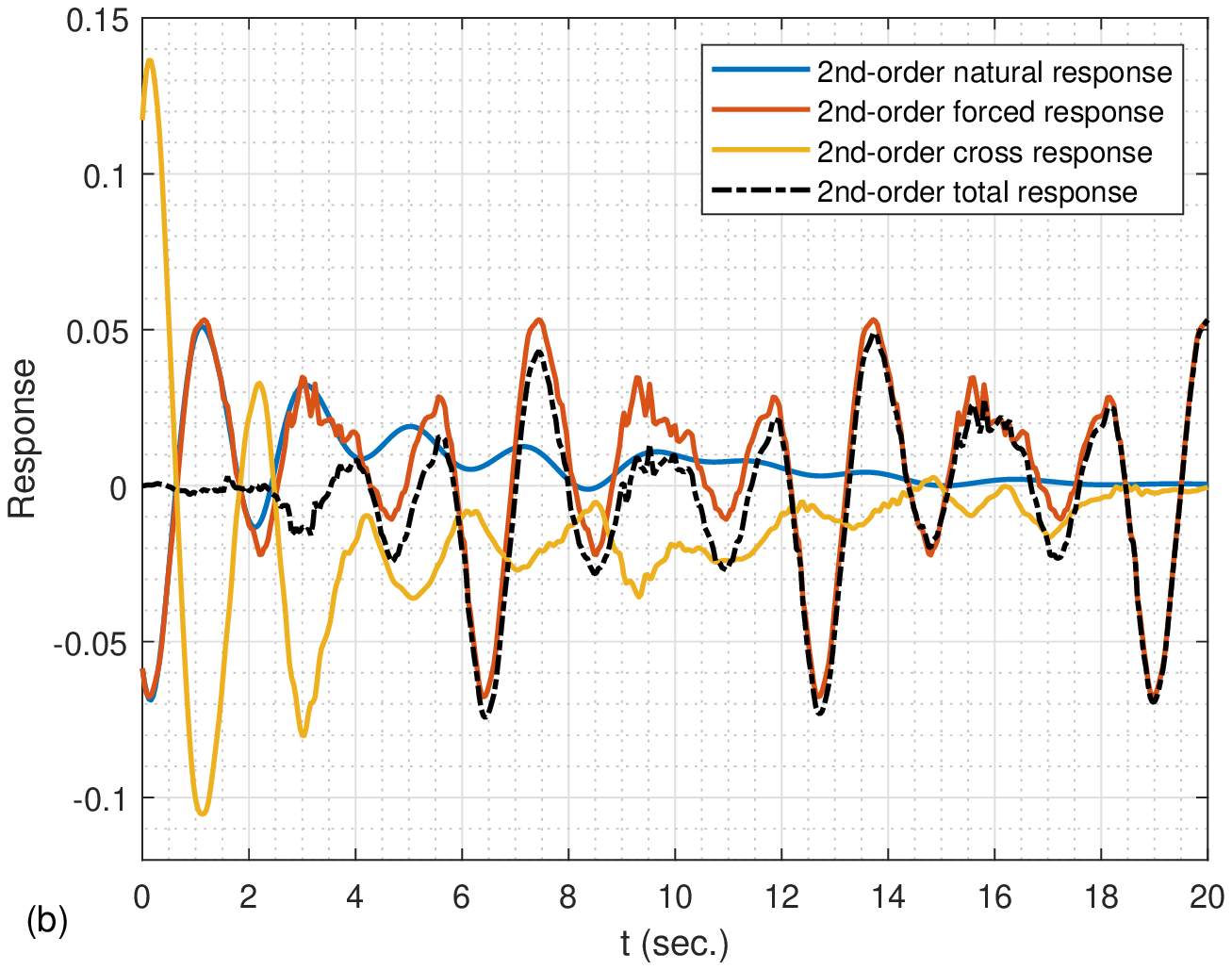}
\end{minipage}
\centering
\caption{Response components (a) for the first-order response and (b) for the second-order response}\label{Exam2_output3_component}
\end{figure}

A regular excitation $f(t) = \sin(\pi t)$ and an irregular excitation $f(t) = \sum_{n=1}^{N_f} A_n \cos(\Omega_n t + \theta_n)$ with $A_n=0.3$ and $\Omega_n$ varying from 0 to 40 with equal interval $1$ are chosen as input excitations. The predicted responses, along with results obtained by the fourth-order Runge--Kutta method, are shown in Fig.~\ref{Exam2_output23}. In both cases, the proposed method accurately predicts system responses.
As presented in Eq.~\ref{y5}, a nonlinear response is the sum of three terms: natural response $y_s(t)$, forced response $y_f(t)$ and cross response $y_c(t)$. These individual terms, as well as their sum to two excitations, are shown in Figs.~\ref{Exam2_output2_component} and \ref{Exam2_output3_component}, respectively. As shown in Figs.~\ref{Exam2_output2_component} and \ref{Exam2_output3_component}, both first- and second-order responses include the natural response $y_s(t)$ and the forced response $y_f(t)$, but the cross response $y_c(t)$ only exists in second-order responses. When $t$ becomes larger, both $y_s(t)$ and $y_c(t)$ diminish due to the presence of system damping, and the total response is entirely governed by $y_f(t)$. Moreover, we notice some features at $t=0$ for these components, including $y_s(0)=-y_f(0)$ for the first-order response and $y_s(0)+y_f(0)=-y_c(0)$ for the second-order response, which are due to imposed zero initial conditions.


\section{Conclusions}
Considering arbitrary irregular excitations, an efficient generalized pole-residue method to compute the nonlinear dynamic response modelled by the Volterra series was developed. A core of the proposed method was obtaining poles and corresponding coefficients of Volterra kernel functions, then those of each order response modelled by each order Volterra series. Once the poles and corresponding coefficients of Volterra kernel functions and excitations were both available, the remaining derivation could follow a similar pole-residue method that had been developed for ordinary linear oscillators. To obtain the poles and corresponding coefficients of Volterra kernel functions, two steps were included: (1) using Laguerre polynomials to decouple higher-order Volterra kernel functions with respect to time and (2) obtaining poles and corresponding coefficients of Laguerre polynomials in the Laplace domain. Because the proposed method gave an explicit, continuous response function of time, it was much more efficient than traditional numerical methods. Moreover, many meaningful physical and mathematical insights were gained because not only each order response but also the natural response, the forced response and the cross response of each order were obtained in the solution procedure. To demonstrate that the proposed method was not only suitable for a system with a known equation of motion but also applicable to a system with an unknown equation of motion, two numerical studies were conducted. For each study, regular excitations and complex irregular excitations with different parameters were investigated. The efficiency of the proposed method was verified by the fourth-order Runge--Kutta method. This paper only computes the response under zero initial conditions. The response under non-zero initial conditions will be investigated in our future work.

\section*{Acknowledgements}
The research was financially supported by the National Natural Science Foundation of China (Grant Nos. 52101302 and 52101339) and the Postdoctoral Research Foundation of China (Grant No. 2021M690521).

\section*{Code availability}
All code that support the findings of this study are available from the authors by email.

%

\section*{References}
\bibliographystyle{unsrt}
\bibliography{refer}

\end{document}